\documentclass[11pt]{amsart}

\ifx\shlhetal\undefinedcontrolsequence\let\shlhetal\relax\fi

\usepackage{amsmath} 
\usepackage{amssymb}

\newtheorem{theorem}{Theorem}[section]

\newtheorem{lemma}[theorem]{Lemma} 
\newtheorem{proposition}[theorem]{Proposition} 
 
\newtheorem{observation}[theorem]{Observation}
\newtheorem{mic}[theorem]{The Main Isomorphism Theorem} 

\theoremstyle{definition}
\newtheorem{definition}[theorem]{Definition}

\newtheorem{problem}{Problem}

\theoremstyle{remark}
\newtheorem{remark}[theorem]{Remark}

\newtheorem{conclusion}[theorem]{Conclusion}

\newtheorem{hypothesis}[theorem]{Hypothesis}
\newtheorem{notation}[theorem]{Notation}

\newcommand{\bA}{{\bf A}}
\newcommand{\cD}{{\mathcal D}}
\newcommand{\bbF}{{\mathbb F}}
\newcommand{\cF}{{\mathcal F}}
\newcommand{\bG}{{\bf G}}
\newcommand{\cK}{{\mathcal K}}
\newcommand{\bL}{{\bf L}}
\newcommand{\cL}{{\mathcal L}}
\newcommand{\bm}{{\bf m}}
\newcommand{\cT}{{\mathcal T}}
\newcommand{\gp}{{\mathfrak p}}
\newcommand{\bV}{{\bf V}}
\newcommand{\bbZ}{{\mathbb Z}}

\newcommand{\coh}{{{\mathbb C}_{\aleph_3}}}
\newcommand{\coA}{{{\mathbb C}_{\bA}}}
\newcommand{\coAr}{{{\mathbb C}_{\bA^r}}}
\newcommand{\coAb}{{{\mathbb C}_{\bA\cap\beta}}}
\newcommand{\cf}{{\rm cf}}
\newcommand{\App}{{{\mathbb A}{\rm pp}}}
\newcommand{\forces}{\Vdash} 
\newcommand{\Dom}{{\rm Dom}}
\newcommand{\Rang}{{\rm Rang}}
\newcommand{\aprod}{\mathop{{\prod}^\bA}}
\newcommand{\aqprod}{\mathop{{\prod}^{\bA^q}}}
\newcommand{\arprod}{\mathop{{\prod}^{\bA^r}}}
\newcommand{\abprod}{\mathop{{\prod}^{\bA\cap\beta}}}
\newcommand{\xtheta}{\vartheta}
\newcommand{\ytheta}{\name{\vartheta}^*}
\newcommand{\ztheta}{\varrho}

\newcount\skewfactor
\def\mathunderaccent#1#2 {\let\theaccent#1\skewfactor#2
\mathpalette\putaccentunder}
\def\putaccentunder#1#2{\oalign{$#1#2$\crcr\hidewidth
\vbox to.2ex{\hbox{$#1\skew\skewfactor\theaccent{}$}\vss}\hidewidth}}
\def\name{\mathunderaccent\tilde-3 }

\title[Viva III]{Vive la Diff\'erence III}
\author{Saharon Shelah}
\address{Institute of Mathematics\\
 The Hebrew University of Jerusalem\\
 Jerusalem 91904, Israel\\
 and  Department of Mathematics\\
 Rutgers University\\
 New Brunswick, NJ 08854, USA}
\email{shelah@math.huji.ac.il}
\urladdr{http://www.math.rutgers.edu/$\sim$shelah}
\thanks{This research was partially supported by the United
States-Israel Binational Science Foundation. Publication 509}  

\subjclass{}
\keywords{Forcing, ultrapowers, strong independence property, bigness
notions, definability} 

\begin{document}

\begin{abstract}
We show that, consistently, there is an ultrafilter $\cF$ on $\omega$ such
that if $N^\ell_n=(P^\ell_n\cup Q^\ell_n,P^\ell_n,Q^\ell_n,R^\ell_n)$ (for
$\ell=1,2$, $n<\omega$), $P^\ell_n\cup Q^\ell_n \subseteq\omega$, and
$\prod\limits_{n< \omega} N^1_n/\cF\equiv\prod\limits_{n<\omega}N^2_n/ 
\cF$  are models of the canonical theory $t^{\rm ind}$ of the strong
independence property, then every isomorphism from $\prod\limits_{n<\omega}
N^1_n/\cF$ onto $\prod\limits_{n<\omega} N^2_n/\cF$ is a product isomorphism.
\end{abstract}

\maketitle

\setcounter{section}{-1}
\section{Introduction} 
In a previous paper \cite{Sh:326} we gave two constructions of models of set 
theory in which the following isomorphism principle fails in various strong
respects:
\begin{description}
\item[(Iso 1)] If $M$, $N$ are countable elementarily equivalent
structures and $\cF$ is a non-principal ultrafilter on $\omega$, then the
ultrapowers $M^*$, $N^*$ of $M$, $N$ with respect to $\cF$ are
isomorphic.
\end{description}
As is well known, this principle is a consequence of the Continuum
Hypothesis.  
Recall that Keisler celebrated theorem ( from \cite{Ke67}) says  that two
models, $M,N$ of cardinality at most $\lambda^+$  are elementarily
equivalent iff for some ultrafilter $\cF$ of $\lambda$, the ultrapowers  
$M^\lambda \cF , N^\lambda / \cF$ are isoorphic, this has given
an algebraic characterization of elementary equivalence.

In \cite{Sh:405} our aim originally was to give a related
example in connection with the well-known isomorphism theorem of Ax and
Kochen. In its general formulation, that result states that a fairly broad
class of Henselian fields of characteristic zero satisfying a completeness
(or saturation) condition are classified up to isomorphism by the structure
of their residue fields and their value groups. The case that interest us
in the 
second paper in this series \cite{Sh:405}, was:  
\begin{description}
\item[(Iso 2)]
If $\cF$ is a non-principal ultrafilter on $\omega$, then the ultraproducts
$\prod\limits_p \bbZ_p/\cF$ and $\prod\limits_p\bbF_p[[t]]/\cF$ are
isomorphic. 
\end{description}

And more generally:

\begin{theorem}
[See \cite{Sh:405}]
\label{PropA}
It is consistent with the axioms of set theory that there is a non-principal
ultrafilter $\cF$ on $\omega$ such that for any two sequences of discrete
rank 1 valuation rings $(R^i_n)_{n=1,2,\ldots}$ ($i=1,2$) having countable
residue fields, any isomorphism $F:\prod\limits_n R^1_n/\cF\longrightarrow 
\prod\limits_n R^2_n/\cF$ is an ultraproduct of isomorphisms $F_n:R^1_n
\longrightarrow R^2_n$ (for a set of $n$'s contained in $\cF$). In
particular, $\cF$--majority of the pairs $R^1_n$, $R^2_n$ are isomorphic.
\end{theorem}

In the case of the rings $\bbF_p[[t]]$ and $\bbZ_p$, we see that (Iso 2)
fails. For this our main work was to show the following statement which
actually from model theoretic point of view is more basic and interesting. 

\begin{theorem}
[See \cite{Sh:405}]
\label{PropB}
It is consistent with the axioms of set theory that there is a nonprincipal
ultrafilter $\cF$ on $\omega$ such that for any two sequences of countable
trees $(T^i_n)_{n=1,2,\ldots}$ for $i=1,2$, with each tree $T^i_n$ countable
with $\omega$ levels, and with each node having at least two immediate
successors, if $\cT^i=\prod\limits_n T^i_n/\cF$, then for any isomorphism
$F:\cT^1\stackrel{\simeq}{\longrightarrow}\cT^2$ there is an element
$a\in\cT^1$ such that the restriction of $F$ to the cone above $a$ is the
restriction of an ultraproduct of maps $F_n: T^1_n\longrightarrow T_n^2$. 
\end{theorem}

\relax From a model theoretic point of view this still is not the right level of
generality for a problem of this type. There are two natural ways to pose
the problem: 

\begin{problem}
\label{Problem1}
Characterize the pairs of countable models $M$, $N$ such that in some
forcing extension, $\prod\limits_n M/\cF\not\simeq\prod\limits_n N/\cF$
for some non-principal ultrafilter $\cF$. 
\end{problem}

\begin{problem}
\label{Problem2}
Characterize the pairs of countable models $M$, $N$ with non-isomorphic
ultrapowers $\mod\omega$ in some forcing extension. (I.e., such that there
is no forcing extension in which for some non-principal ultrafilter $\cF$ on
$\omega$ we have $M^\omega/\cF\simeq N^\omega/\cF$.)
\end{problem}

\noindent [There are two variants of the second problem: the ultrapowers may
be formed either using one ultrafilter twice (called \ref{Problem2}(A)), or
may consider using any two ultrafilters (called \ref{Problem2}(B)), but see
below.]  

\begin{problem}
\label{Problem3}
Let us write $M\leq N$ whenever in every forcing extension, if $\cF$ is
an ultrafilter on $\omega$ such that $N^\omega/\cF$ is saturated, then
$M^\omega/\cF$ is also saturated. Characterize this relation.
\end{problem}

This is somewhat like the Keisler order (see Keisler \cite{Ke67}, or
\cite{Sh:a}, or \cite[Chapter VI]{Sh:c}), but does not depend on the fact
that the ultrafilter is regular, so some of the results there apply to
Problem \ref{Problem3}, this in turn implies results on Problem
\ref{Problem2}(A). We can replace $\aleph_0$ here by any cardinal $\kappa$
satisfying $\kappa^{<\kappa}=\kappa$. 

Now, by \cite{Sh:13}, there is an ultrafilter $\cD$ on $2^{\aleph_0}$ such
that for ountable models $M,N$ 
\[M\equiv N\quad\Rightarrow\quad M^{2^{\aleph_0}}/\cD\simeq N^{2^{\aleph_0}}
/\cD.\]
Also, if $2^{\aleph_0}=\aleph_1$, $\cF$ is a non-principal ultrafilter on
$\omega$ and $M_1\equiv M_2$ are countable, then $M^\omega_1/\cF\simeq
M^\omega_2/\cF$ (as they are saturated); similarly if $M^\ell_n$ are
countable models (for $\ell=1,2$, $n<\omega$), $M_\ell=\prod\limits_{n<
\omega} M^\ell_n/\cF_\ell$, and $\cF_\ell$ are non-principal ultrafilters on
$\omega$, then $M_1\equiv M_2\ \Rightarrow\ M_1\cong M_2$. On the other
hand, if $2^{\aleph_0}>\aleph_1$, then by \cite[Ch VI]{Sh:c} for every
regular cardinal $\theta$, $\aleph_1\leq\theta<2^{\aleph_0}$ we have an
ultrafilter $\cF_\theta$ on $\omega$ such that the down cofinality of
$(\omega,<)^\omega/\cF_\theta$ above $\omega$ is $\theta$ (so $\theta_1\neq
\theta_2\ \Rightarrow\ (\omega,<)^\omega/\cF_{\theta_1}\not\simeq
(\omega,<)^\omega/\cF_{\theta_2}$).  
 
The present paper is dedicated to
hadding some
further
 light. Working on \cite{Sh:405}
we had hoped to continue it sometime. However, we actually began only when 
Jarden asked:
\begin{enumerate}
\item[$(*)$] Suppose that $F_n^\ell$ are finite fields (for $n<\omega$,
$\ell=1,2$). Can we have (a universe and) an ultrafilter $\cF$ on $\omega$
such that $\prod\limits_{n<\omega} F_n^1 /\cF$
and $\prod\limits_{n<\omega} 
F_n^2/\cF$?  are elementaily equivalent but no isomorphic.
\end{enumerate}

That was not an arbitrary question: he knew that many such 
pairs of ultraproducts are elementaily equivalent,
because the first order theory of a field $F$
which is isomorphic to an ultraproduct  of finite fields
is determine by its chracteristic  and its
subfield of algebraic elements. Hencewe can find an equivalence
relation $E_k$  on the family of finite fields
for $k<w$ ,
each  with finitely many equivalence classes
such that if $ F^1_n,F^2_n$ are finite fields for $n<\omega$
and $\cF$ is a non pricipal ultrafilter on $\omega$
and for each k the set $\{n<\omega : (F^1_n)  E_k (F^2_n )
\} $ belongs to $\cF$ then the
respective  ultraproducts
are isomorphic.

Jarden asked me, I inquire whether it has the strong independence property
and told him what it is, he says yes. Years later finishing the work
on the paper he deny any knowledge on this, and this is my recollection.
Cherlin, to whom he refer me, give me the 
reference to the strong independence property for finite field:
Duret \cite[pp. 136--157]{Du80}.

Here we continue \cite[\S 3]{Sh:326}, \cite[\S 1]{Sh:405}. To give an
affirmative answer to $(*)$, we show that after adding $\aleph_3$ Cohen
reals to a suitable ground model, one gets a universe with an ultrafilter 
$\cF$ on $\omega$ and a sequence of models $\langle M_n:n<\omega
\rangle$ on $\omega$ such that  
\begin{enumerate}
\item[$(**)$] \underline{if} $N^\ell_n=(P^\ell_n\cup Q^\ell_n,P^\ell_n,
Q^\ell_n,R^\ell_n)$ (for $\ell=1,2$, $n<\omega$), $P^\ell_n\cup Q^\ell_n
\subseteq\omega$, and $\prod\limits_{n<\omega} N^1_n/\cF\equiv\prod\limits_{
n<\omega}N^2_n/\cF$ are models of the canonical theory $t^{\rm ind}$ of the
strong independence property (see Definition \ref{bigind}), 

\noindent\underline{then} every isomorphism from $\prod\limits_{n<\omega}
N^1_n/\cF$ onto $\prod\limits_{n<\omega} N^2_n/\cF$ is (first order)
definable in $\prod\limits_{n<\omega} M_n/\cF$
for some expansions $M_n$ of  $N^1_n,N^2_n$ simultaneously,
or what is equivalent buthopefully more transparent 
if $F$ is an isomorphism from  
 $N^1 = \prod\limits_{n<\omega}
N^1_n/\cF$
onto
 $N^2= \prod\limits_{n<\omega}
N^3_n/\cF$  then we can find unary functions $F_n$
from $N^1_n$ onto $N^2_n$ for every $n<\omega $
such that the set of $n$ for which $F_n$ is an isomorphism
from$N^1_n$ onto $N^2_n$
belongs to the ultrafilter
and  $\prod\limits_{n<\omega}
(N^1_n , N^2_n , F_n)/\cF$
is $(N^1, N^2,F)$
 . 
\end{enumerate} 

Out forcing is adding $\aleph_3$ Cohen reals, but we 
 need that   our model of set theory, i.e.
 the universe, satisfies some
conditions over which we force.
There are two ways to get a ``suitable'' ground model. The first way
involves taking any ground model which satisfies a portion of the GCH, and
extending it by an appropriate preliminary forcing, which generically adds
the {\em name} for an ultrafilter which will appear after addition of the
Cohen reals. The alternative approach, which we consider
more
 model--theoretic,
is to start with an $\bL$--like ground model and use instances of diamond
(or related weaker principles) to prove that a sufficiently generic name
already exists in the ground model. We will fully present the first approach
- the second one should be then an easy modification of the arguments
presented in \cite[\S 1]{Sh:405}. 

Our presentation is slightly more general than needed for (**). By allowing
more 
what we call
"bigness"
 properties to be involved in the definition of $\App$, we leave
room for getting analogs of (**) for more classes of models (getting the
conclusion for all of them at once, or possibly only for some) - as long as
the respective bigness notions are like in \ref{verygoodb}. This, we hope, 
would be helpful in connection with Problems \ref{Problem1}, \ref{Problem2}
above. For the problem on fields only the case associated with the strong
independence property is needed; 
general bigness notions appear for possible
general treatment.

Let us comment on our general point of view. In this paper we try to advance
in Problems \ref{Problem1}+\ref{Problem2}(A) and for this, it seemed, we can
take the maximal $\name{\Gamma}$, i.e., allow all $\aleph_0$-bigness
notions.
 However, concerning Problem \ref{Problem3} (investigating the
partial order $\leq$ or models), for showing $M\nleq  N$, the
construction causes $N^\omega/\cF$ to be almost always non
$\aleph_3$--saturated. We need stronger tools for them.  
\medskip

The two previous papers benefited from Gregory Cherlin, the present one
benefited from Andrzej Ros{\l}anowski, thank you!
\medskip

We continue this investigations in \cite{Sh:F503}.

\begin{notation}
Our notation is standard and compatible with that of classical textbooks
(like Chang and Keisler \cite{CK} and Jech \cite{J}). In forcing we keep the
older convention that {\em a stronger condition is the larger one}.   
\begin{enumerate}
\item We will use two forcing notions denoted by $\coh$ and $\App$ (see
Definitions \ref{generalia} and \ref{app}, respectively). Conditions in
these forcing notions will be called $p,q,r$ (with possible
sub/super-scripts). 
\item All names for objects in forcing extensions will be denoted with a
tilde below (e.g., $\name{a}$, $\name{\gp}$).
\item The letter $\tau$ (with possible sub/super-scripts) stands for a
vocabulary of a first order language; me may also write $\tau(M)$,
$\tau(T)$ for a model $M$ or theory $T$ with the obvious meaning. We will
use letter $\gp$ (with sub/super-scripts) to denote types.  
\item The universe of a model $M$ will be denoted $|M|$, but we will often
abuse this notation and write, e.g., $a\in M$. The cardinality of a set $A$
will be denoted $\|A\|$, and, for a model $M$, $\|M\|$ will stand for the
cardinality of its universe.  
\end{enumerate}
\end{notation}

\section{Bigness notions}
In this section we will quote relevant definitions and results from
\cite[Chapters X, XI]{Sh:e} (=\cite{Sh:384}, \cite{Sh:482}), but we somewhat
restrict ourselves here. The reader interested in the field case only may
jump directly to Definition \ref{bigind}.

\begin{definition}
[See {\cite[Chapter XI, \S1]{Sh:e}}]
\label{482def1:10}
Let $T$ be a complete first order theory (in a vocabulary $\tau$), and $\cK=
\cK_T$ be a class of models of $T$ partially ordered by $\prec$. Also let
$t$ be a first order theory with a countable vocabulary $\tau(t)$ (including
equality, treating function symbols as predicates).   
\begin{enumerate}
\item We  say that $\cK'$ is an $A$--place in $\cK$ if
\begin{enumerate}
\item[(a)] $\cK'\subseteq\cK$,
\item[(b)] if $M\in\cK'$, then $A\subseteq M$,
\item[(c)] if $M\prec N$ are from $\cK$ and $A\subseteq M$, then $M\in\cK'\
\Leftrightarrow\ N\in\cK'$, 
\item[(d)]if $M\in\cK'$ and $A\subseteq N\in\cK$ and $M,N$ are isomorphic
over $A$, then $M\in\cK'\ \Leftrightarrow\ N\in\cK'$.
\end{enumerate}
\item For $A\subseteq M\in\cK$ we let $\cK'=\cK_{A,M}$ be the class
\[\{N: A\subseteq N\mbox{ and }\bar{a}\in {}^{\omega>}A\ \Rightarrow\ {\rm
tp}(\bar{a},\emptyset,M)={\rm tp}(\bar{a},\emptyset,N)\ \}.\]
We call it {\em the $(A,M)$--place}.
\item {\em A local bigness notion $\Gamma$ for $\cK$\/} (without parameters,
in one variable $x$) is a function with domain $\cK$ which for every model
$M\in\cK$ gives 
\[\begin{array}{l}
\Gamma^-_M=\Gamma^-(M)\subseteq \{\varphi(x,\bar{a}):\varphi\in\cL(\tau)\
\&\ \bar{a}\subseteq M\},\\ 
\Gamma^+_M=\Gamma^+(M)=\{\varphi(x,\bar{a}):\varphi\in\cL(\tau)\ \&\
\bar{a}\subseteq M\}\setminus\Gamma^-_M 
  \end{array}\]
such that
\begin{enumerate}
\item[(a)] $\Gamma^-_M$ is preserved by automorphisms of $M$,
\item[(b)] $\Gamma^-_M$ is a proper ideal, i.e., $\Gamma^+_M\neq\emptyset$
and 
\begin{enumerate}
\item[$(\alpha)$] if $M\models(\forall x)(\varphi(x,\bar{a})\Rightarrow \psi
(x,\bar{b}))$ and $\psi(x,\bar{b})\in\Gamma^-_M$, then
$\varphi(x,\bar{a})\in\Gamma^-_M$, 
\item[$(\beta)$] if $\varphi_1(x,\bar{a}_1),\varphi_2(x,\bar{a}_2)
\in \Gamma^-_M$, then $\varphi_1(x,\bar{a}_1)\vee \varphi(x,\bar{a}_2)\in
\Gamma^-_M$. 
\end{enumerate}
\end{enumerate}
Elements of $\Gamma^-_M$ are called {\em $\Gamma$--small in $M$}, members of 
$\Gamma^+_M$ are {\em $\Gamma$--big}.

{\em A local bigness notion $\Gamma$ for $\cK$ with
parameters\footnote{Alternatively use the monster model.} from $A$} is
defined similarly but $\Dom(\Gamma)$ is an $A$--place $\cK'$ in $\cK$.
\item We say that a local bigness notion $\Gamma$ is {\em invariant\/} if
for $M\prec N$ from $\cK$ we have $\Gamma^-_M\subseteq\Gamma^-_N$ and
$\Gamma^+_M \subseteq \Gamma^+_N$. 
\item {\em A $\Gamma$--big type $\gp(x)$ in $M$} is a set of formulas
$\psi(x,\bar{a})$
all of  whose
 finite conjunctions are $\Gamma$--big in $M$. 
\item A {\em  pre $t$--bigness notion scheme} $\Gamma$ is a sentence
$\psi_\Gamma$ (in possibly infinitary logic) in the vocabulary
$\tau(t)\cup\{P^*\}$, where $P^*$ is a unary predicate. 
\item {\em An interpretation with parameters of $t$ in a model $M\in\cK$} 
is $\bar{\varphi}=\langle\varphi_R(\bar{y}_R,\bar{a}_R):R\in\tau(t)\rangle$,
where $\varphi_R\in\cL(\tau)$ and $\bar{a}$ is a sequence of appropriate
length of elements of $M$. So $R$ is interpreted as  
\[\{\bar{b}:M\models\varphi_R(\bar{b},\bar{a}_R),\ \lg(\bar{b})=
\lg(\bar{y}_R)\ (=\mbox{ the arity of }R)\;\}.\]
The interpreted model is called $M[\bar{\varphi}]$ and we demand that it is
a model of $t$. 
\item For a pre $t$--bigness notion scheme $\Gamma=\psi_\Gamma$ and an
interpretation $\bar{\varphi}$ of $t$ in $M\in\cK$ with parameters from
$A\subseteq M$, we define the $\bar{\varphi}$--derived {\em local bigness
notion $\Gamma[\bar\varphi]$ with parameters from $A\subseteq M$ (in the
$\cK_{A,M}$--place)} as follows:\\  
Given $M'\in\cK_{(A,M)}$. A formula $\vartheta(x,\bar{b})$ in $\cL(\tau)$
(with parameters from $M$, of course) is $\Gamma[\bar{\varphi}]$--big in $M$
if for any quite saturated $N^*$, $M\prec N^*$, letting  
\[P^*=\{a\in N^*[\bar{\varphi}]:N^*\models\vartheta[a,\bar{b}]\}\]
we have $(N^*[\bar{\varphi}],P^*)\models\psi_\Gamma$.  

We write $\Gamma=\Gamma_{(\bar{\varphi},t,\psi_\Gamma)}$. 

\item We omit the ``pre'' if every $\Gamma[\bar{\varphi}]$ is an invariant
local bigness notion (for our fixed $\cK$). So it is enough in (8) above if
we define $\Gamma_M$ when $M\prec M'$.
\end{enumerate}
\end{definition}

\begin{proposition}
\label{1.1A}
If $\Gamma$ is a local bigness notion for $\cK$ with parameters in $A$,
$M\in\cK_{A,M'}$ and $\gp(x)$ is a $\Gamma$--big type in $M$,
\underline{then} it can be extended to $\Gamma$--big notion ${\mathfrak q}$
in $M$ which is a complete type over $M$.
\end{proposition}

\begin{proposition}
For $T,\cK=\cK_T$ and $t$ as in \ref{482def1:10},
\begin{enumerate}
\item[$(\boxtimes)$] if $N\prec M$ are from $\cK$, and $\bar{\varphi} 
=\langle\varphi_R:(\bar{y}_R,\bar{a}_R):R\in\tau(t)\rangle$ is an
interpretation of $t$ in $N$, 
\underline{then} $\bar{\varphi}$ is an interpretation of $t$ in $M$
(i.e., $M[\bar{\varphi}]\models t$). 
\end{enumerate}
\end{proposition}

The following definition illuminates the most important aspect of Definition
\ref{bigind} (which is central for our present paper), and also it is needed
for more general results.

\begin{definition}
\label{verygoodb}
Let $t$ be a first order theory in a vocabulary $\tau(t)$. Suppose that
$\Gamma$ is a $t$--bigness notion scheme, $P\in\tau(t)$ is a unary 
predicate, and $\xtheta(y,x)$ is a $\tau(t)$--formula. We say that $\Gamma$
is {\em $(\aleph_2,\aleph_1)$--$(P,\xtheta)$--separative with a witness $X$}
whenever the following condition $(\circledast)^{P,\xtheta}_\Gamma$ holds.  
\begin{enumerate}
\item[$(\circledast)^{P,\xtheta}_\Gamma$] Assume that $M$ is an
$\aleph_2$--compact $\tau$--model, $\bar{\varphi}=\langle\varphi_R:
(\bar{y}_R,\bar{a}_R):R\in\tau(t)\rangle$ is an interpretation of $t$ in
$M$. Then $X\subseteq |M|$ is of cardinality at most $\aleph_1$, includes
all parameters of $\bar{\varphi}$ and 
\begin{enumerate}
\item[{\bf if}] $N\prec M$, $X\subseteq |N|$, $\|N\|\leq\aleph_1$, and
$\gp(x)$ is a $\Gamma[\bar{\varphi}]$--big type over $N$, $\|\gp(x)\|\leq
\aleph_1$, and $a_1,a_2$ are distinct members of $|M|\setminus |N|$ with
\[M\models \varphi_P[a_1]\ \wedge \varphi_P[a_2]\]
\item[{\bf then}] the type $\gp(x)\cup\{\xtheta(a_1,x)\equiv\neg\xtheta(a_2,
x)\}$ is $\Gamma[\bar{\varphi}]$--big.
\end{enumerate}
We may omit the ``$(\aleph_2,\aleph_1)$--''.
\end{enumerate}
\end{definition}

\begin{definition}
[See {\cite[Def. 3.4, 3.5, Chapter XI]{Sh:e}}]
\label{bigind}
\begin{enumerate}
\item $t^{\rm ind}=t^{\rm ind}_0$ is the first order theory in vocabulary
$\tau(t^{\rm ind})=\{P,Q,R\}$, where $P,Q$ are unary predicates and $R$ is a
binary predicate, including formulas  
\[\begin{array}{l}
(\forall x)(\forall y)(x\;R\;y\ \Rightarrow\ P(x)\ \wedge\ Q(y)),
\quad\mbox{ and}\\
(\forall x)(P(x)\vee Q(x))
  \end{array}\]
and saying that for each $n<\omega$ and any distinct elements $a_1,\ldots,
a_{2n}\in P$, there is $c\in Q$ such that
\[a_i\;R^M\; c\quad\mbox{if and only if}\quad i\leq n.\]
$t^{\rm ind}_1$ is $t^{\rm ind}$ plus
\[(\forall x)(\forall y)(\exists z)\big(Q(x)\ \wedge\ Q(y)\ \wedge\ x\neq y\
\Rightarrow\ P(z)\wedge (z\;R\;x\equiv \neg z\;R\;y)\big).\]
\item We define a pre $t^{\rm ind}$--bigness notion scheme $\Gamma^{\rm
ind}$ as follows. The sentence $\psi_{\Gamma^{\rm ind}}$ says that
$P^*\subseteq Q$ and $(P,Q,R,P^*)$ satisfies:
\begin{quotation}
for every $n<\omega$, there is a finite set $A\subseteq P$ such that for
every distinct $a_1,\ldots,a_{2n}\in P\setminus A$ there is $c\in P^*$
satisfying 
\[a_\ell\;R\;c\quad \mbox{for }\ell\leq n,\quad\mbox{ and }\quad \neg
a_\ell\;R\;c \quad \mbox{for }n<\ell\leq 2n.\]
\end{quotation}
(So $\psi_{\Gamma^{\rm ind}}$ is not first order.)
\item We say that a first order theory $T$ has the strong independence
property if some formula $\xtheta(x,y)$ define a two place relation
which is a model of $t^{\rm ind}_1$  with $P,Q$ choosen
as $x+x$
\end{enumerate}
\end{definition}

Plainly, 

\begin{proposition}
\label{1.3A}
\begin{enumerate}
\item For a model $M$ of $t^{\rm ind}_1$, an automorphism $\pi$ of $M$ is
determined by $\pi\restriction P^M$ (i.e., if $\pi_1,\pi_2\in{\rm Aut}(M)$
are such that $\pi_1\restriction P^M=\pi_2\restriction P^M$, then $\pi_1=
\pi_2$).
\item Moreover, if $\bar{\varphi}$ is an interpretation of $t^{\rm ind}_1$
in $M^*$, $M=M^*[\bar{\varphi}]$, $\pi\in {\rm Aut}(M)$ and $\pi\restriction
P^M$ is definable with parameters in $M^*$, \underline{then} so is $\pi$. 
\end{enumerate}
\end{proposition}

\begin{proposition}
[See {\cite[Chapter XI, \S3]{Sh:e}}]
$\Gamma^{\rm ind}$ is a $t^{\rm ind}$--bigness notion scheme. It is
$(\aleph_2,\aleph_1)$--$(P,\xtheta)$--separative with witness $\emptyset$,
where $P$ is given and $\xtheta(y,x)=y\;R\;x$.   
\end{proposition}

\begin{definition}
\label{1.5}
A mapping $F:N^1\longrightarrow N^2$ is {\em a $\Delta$--embedding from
$N^1$ to $N^2$} whenever $\Delta$ is a set of formulas in $L_{\omega,\omega}
(\tau(N^1)\cap \tau(N^2))$ such that
\begin{enumerate}
\item[{\bf if}]  $\varphi\in\Delta$ and $N^1\models\varphi[a_1,\ldots,a_n]$,
\item[{\bf then}] $N^2\models\varphi[F(a_1),\ldots,F(a_n)]$.
\end{enumerate}
[If $\Delta$ is closed under negation, then we have ``if and only if''.]
\end{definition}

\section{The forcing notion $\App$}
As explained in the introduction, we work in a Cohen generic extension of a
suitable ground model. In this section we present how that ground model can
be obtained: we start with $\bV\models{\rm GCH}$ and we force with the
forcing notion $\App$ defined in \ref{app} below
, the App comes for approximations, as themebers are approxiamtions
to a  name for an  ultrafilter as we desire.

\begin{definition}
\label{generalia}
\begin{enumerate}
\item The Cohen forcing adding $\aleph_3$ Cohen reals is denoted by
$\coh$. Thus a condition $p$ in $\coh$ is a finite partial function from
$\aleph_3\times\omega$ to $\omega$, and the order of $\coh$ is the natural
one. The canonical $\coh$--name for $\beta^{\rm th}$ Cohen real will be
called $\name{x}_\beta$. 
\item Let $\bA\subseteq\aleph_3$. For a condition $p\in\coh$, its
restriction to $\bA\times\omega$ is called $p\restriction\bA$, and we let
$\coh\restriction\bA=\coA=\{p\restriction\bA:p\in\coh\}$. Also, we let
$\name{\omega}^*_{\bA}=({}^{\textstyle\omega}\omega)^{\bV^{\coh\restriction
\bA}}$. 
\item For a sequence $\langle A_n:n<\omega\rangle$ of non-empty sets (and
$\bA\subseteq\aleph_3$), we define  
\[\begin{array}{r}
\aprod\limits_{n<\omega} A_n=\{f\in\bV^{\coh\restriction\bA}:f\mbox{ is a
function with domain }\omega,\ \ \\
\mbox{and such that }f(n)\in A_n\mbox{ for all }n\;\},
  \end{array}\] 
and similarly for models. 
\item For $\bA\subseteq\aleph_3$ and $m<\omega$, let $I^m_{\bA}$ be the set
of all $\omega$--sequences of canonical $\coA$--names for subsets of
$\omega^m$. Let $Q_{\bar{s}}$ (for $\bar{s}\in I^m_{\bA}$, $m<\omega$) be an
$m$-ary predicate, $Q_{\bar{s}_0}\neq Q_{\bar{s}_1}$ whenever $\bar{s}_0\neq
\bar{s}_1$, and let 
\[\tau_{\bA}=\{Q_{\bar{s}}:\bar{s}\in I^m_{\bA}\ \&\ m<\omega\}\]
(so $\|\tau_{\bA}\|=\aleph_1\cdot \|\bA\|$). Let $\name{M}^n_{\bA}$ be a
$\coA$--name for the $\tau_{\bA}$--model on $\omega$ such that if $\bar{s}= 
\langle \name{s}_n:n<\omega\rangle\in I^m_{\bA}$, then $\forces_{\coA}
(Q_{\bar{s}})^{\name{M}^n_A}=\name{s}_n$. 
\end{enumerate}
\end{definition}

\begin{definition}
\label{2.1A}
\begin{enumerate}
\item A function $\bG$ is called {\em an $(\aleph_3,\aleph_2)$--bigness
guide} if the domain $\Dom(\bG)$ of $\bG$ is
\[\begin{array}{ll}
\{(\bA,\name{\cF}):&\bA\subseteq\aleph_3,\ \|\bA\|\leq\aleph_1,\mbox{ and}\\
&\name{\cF}\mbox{ is a $\coA$--name of a non principal ultrafilter on }
\omega\;\},
  \end{array}\]
and 
\begin{enumerate}
\item[$(\alpha)$] $\bG(\bA,\name{\cF})$ is a set of triples $(t,
\name{\Gamma},\name{\bar{\varphi}})$, where $t$ is a first order theory (or
just a $\coA$--name of a first order theory), $\name{\Gamma}$ is a
$\coA$--name of $t$--bigness notion scheme, and $\name{\bar{\varphi}}$ is (a
$\coA$--name for) an interpretation of $t$ in $\aprod\limits_{n<\omega}
\name{M}_\bA^n/\name{\cF}$, and $\|\bG(\bA,\name{\cF})\|\leq\aleph_2$, and 
\item[$(\beta)$] if $(\bA^\ell,\name{\cF}^\ell)\in\Dom(\bG)$ for $\ell=1,2$, 
$\bA^1\subseteq\bA^2$ and $\forces_{{\mathbb C}_{\bA_2}}\name{\cF}^1
\subseteq\name{\cF}^2$, then $\bG(\bA^1,\name{\cF}^1)\subseteq\bG(\bA^2, 
\name{\cF}^2)$.
\end{enumerate}
\item An $(\aleph_3,\aleph_2)$--bigness guide $\bG$ is {\em ind--full\/} if 
\begin{enumerate}
\item[$(\gamma)$] for every $(\bA,\name{\cF})\in\Dom(\bG)$ and a ${\mathbb
C}_\bA$--name $\name{\bar{\varphi}}$ for an interpretation of $t^{\rm ind}$
in $\aprod\limits_{n<\omega}\name{M}_\bA^n/\name{\cF}$ we have $(t^{\rm
ind},\Gamma^{\rm ind},\name{\bar{\varphi}})\in\bG(\bA,\name{\cF})$. 
\end{enumerate}
\item We say that $\bG$ is {\em full} whenever the following condition holds.
\begin{enumerate}
\item[$(\boxplus)$] Assume $(\bA,\name{\cF})\in\Dom(\bG)$ and $\name{t}$ is
a $\coA$--name of a first order theory in the vocabulary $\tau(\name{t})\in
{\mathcal H}(\aleph_1)$, $\name{\psi}$ is a $\coA$--name for a pre
$\name{t}$--bigness notion scheme, $\name{\psi}\in L_{\aleph_1,\aleph_1}(
\tau(\name{t})\cup\{P^*\})$. Let $\name{\bar{\varphi}}$ be a $\coA$--name
for an interpretation of $\name{t}$ in $\aprod\limits_{n<\omega}
\name{M}_\bA^n/\name{\cF}$. Suppose also that $\name{\Delta}$ is (a
$\coA$--name for) a set of $L_{\omega,\omega}(\tau(\name{t}))$--formulas
such that $(\name{\bar{\varphi}},\name{t},\name{\psi})$ defines a bigness
notion $\Gamma=\Gamma_{(\name{\bar{\varphi}},\name{t},
\name{\psi})}$.\footnote{We can fix a ${\mathbb C}_{\aleph_3}$--name of
countable first order theory; really $\name{\cF}$ serves simultaneously for
all.} 
\underline{Then} $(\name{t},\Gamma,\name{\bar{\varphi}})\in\bG(\bA,
\name{\cF})$.
\end{enumerate}
[The main case for us is $t=t^{\rm ind}$, $\Gamma=\Gamma^{\rm ind}$.]
\end{enumerate}
\end{definition}

The clause \ref{2.1A}(2) is added for our particular application. It can be
replaced by use of different bigness notions.   

\begin{proposition}
\label{2.1B}
\begin{enumerate}
\item There is a full $(\aleph_3,\aleph_2)$--bigness guide $\bG$.
\item If a bigness guide $\bG$ is full, \underline{then} it is ind--full.
\end{enumerate}
\end{proposition}

\begin{proof}
Trivial.
\end{proof}

\begin{definition}
\label{app}
Let $\bG$ be an $(\aleph_3,\aleph_1)$--bigness guide. We define the forcing
notion $\App=\App_{\bG}$. (When $\bG$ is fixed, as typically in the present
paper, we do not mention it.)   
\begin{enumerate}
\item {\bf A condition $q$ in $\App$} is a triple $q=(\bA,\name{\cF},
\name{\bar{\Gamma}})=(\bA^q,\name{\cF}^q,\name{\bar{\Gamma}}^q)$ such that:
\begin{enumerate}
\item[(a)] $\bA$ is a subset of $\aleph_3$ of cardinality $\leq\aleph_1$;
\item[(b)] $\name{\cF}$ is a canonical $\coA$--name of a non-principal
ultrafilter on $\omega$, such that for $\beta\in\bA$,  
\[\name{\cF}\restriction(\bA\cap\beta)\stackrel{\rm def}{=}\name{\cF}\cap\{
\name{a}:\name{a}\mbox{ is a $\coAb$--name of a subset of }\omega\;\}\]
is a $\coAb$--name (of an ultrafilter on $\omega$);
\item[(c)] $\name{\bar{\Gamma}}=\langle\name{\Gamma}_\beta:\beta\in\bA\ \&\
\cf(\beta)=\aleph_2\rangle$, where each $\name{\Gamma}_\beta$ is a local
bigness notion $\Gamma[\name{\bar{\varphi}}]$ for some $(t,\Gamma,\name{
\bar{\varphi}})\in\bG(\bA\cap\beta,\name{\cF}\restriction (\bA\cap\beta))$;
\item[(d)] If $\cf(\beta)=\aleph_2$, $\beta\in\bA$, then it is forced (i.e.,
$\forces_{{\mathbb C}_{\aleph_3}}$) that:\\ 
the type realized by $\name{x}_\beta$ over the model
$\abprod\limits_{n<\omega}\name{M}^n_{\bA\cap\beta}/\big(\name{\cF}
\restriction(\bA\cap\beta)\big)$ (so it is a type in the
vocabulary $\tau_{\bA\cap \beta}$) is $\name{\Gamma}_\beta$--big complete,
and moreover this type is a $\coAb$--name. We call it ``the  type induced by
$\name{x}_\beta$ according to $q$''.   
\end{enumerate}
\item {\bf The order ${}\leq_{\App}{}={}\leq{}$ of $\App=\App_\bG$} is the
natural one: \quad $q_1\leq q_2$ \quad if and only if \quad $\bA^{q_1}
\subseteq\bA^{q_2}$, $\forces_{{\mathbb C}_{\bA_2}}\name{\cF}^{q_1}\subseteq 
\name{\cF}^{q_2}$, and $\name{\bar{\Gamma}}^{q_2}\restriction \bA^{q_1}=
\name{\bar{\Gamma}}^{q_1}$.
\item We say that $q_2\in\App$ is {\em an end extension\/} of $q_1\in\App$,
and we write $q_1\le_{\rm end} q_2$, if $q_1\leq q_2$ and $\sup(\bA^{q_1})
\leq\min(\bA^{q_2}\setminus\bA^{q_1})$.
\item For a condition $q\in\App$ and an ordinal $\beta\in\aleph_3$ we define
$q\restriction\beta=\big(\bA^q\cap\beta,\name{\cF}^q\restriction(\bA^q\cap
\beta), \name{\bar{\Gamma}}^q\restriction (\bA^q\cap\beta)\big)$. 
\item For $\beta<\aleph_3$ we let $\App\restriction\beta=\{q\in\App:\bA^q
\subseteq\beta\}$ with inherited order. If $G\subseteq\App$ is generic over
$\bV$, then we let $G\restriction \beta=G\cap(\App\restriction\beta)$. 
\end{enumerate}
\end{definition}

One easily checks that 
\begin{proposition}
\begin{enumerate}
\item If $q\in\App$, $\beta<\aleph_3$, \underline{then} $q\restriction\beta
\in\App$ and $q\restriction\beta\leq_{\rm end} q$.
\item Both $\leq_{\App}$ and $\leq_{\rm end}$ are partial orders on $\App$.
\end{enumerate}
\end{proposition}

\begin{lemma}
\label{1.7Lemma}
If $\langle q_\zeta:\zeta<\xi\rangle$ is an increasing sequence of members
of $\App$, $\xi\leq\aleph_1$, and $q_{\zeta_1}\leq_{\rm end} q_{\zeta_2}$
for $\zeta_1<\zeta_2$, \underline{then} there is $q\in\App$ such that
$\bA^q=\bigcup\limits_{\zeta<\xi}\bA^{q_\zeta}$ and  $q_\zeta\leq_{\rm end}
q$ for all $\zeta<\xi$. 
\end{lemma}

\begin{proof}
We may assume that $\xi>0$ is a limit ordinal. If $\cf(\xi)>\aleph_0$, then
we let $\bA^q=\bigcup\limits_{\zeta<\xi}\bA^{q_\zeta}$, $\name{\cF}^q=
\bigcup\limits_{\zeta<\xi}\name{\cF}^{q_\zeta}$ and $\name{\bar{\Gamma}}^q= 
\bigcup\limits_{\zeta<\xi}\name{\bar{\Gamma}}^{q_\zeta}$. If $\cf(\xi)=
\aleph_0$, then additionally we have to extend $\bigcup\limits_{\zeta<\xi}
\name{\cF}^{q_\zeta}$ to a ${\mathbb C}_{\bA^q}$--name of an ultrafilter on
$\omega$, which is no problem.   
\end{proof}

\begin{lemma}
\label{1.8Lemma}
Suppose that $q\in\App$, $\bA^q\subseteq\gamma\in\aleph_3$, and $\name{\gp}$
is a ${\mathbb C}_{\bA^q}$--name of a type over the model $\aqprod\limits_{n
<\omega}\name{M}^n_{\bA^q}/\name{\cF}^q$ (so in the vocabulary
$\tau_{\bA^q}$, finitely satisfiable in $\aqprod\limits_{n<\omega}
\name{M}^n_{\bA^q}/\name{\cF}^q$). \underline{Then}:
\begin{enumerate}
\item If $\cf(\gamma)<\aleph_2$, then there is a condition $r\in\App$
stronger than $q$ such that $\bA^r=\bA^q\cup\{\gamma\}$, and 
\[\forces_{\coAr}\mbox{`` }\name{x}_\gamma/\name{\cF}^r\mbox{ realizes }
\name{\gp}\mbox{ in }\arprod\limits_{n<\omega}\name{M}^n_{\bA^r}/
\name{\cF}^r\mbox{ ''}.\]
\item If $\cf(\gamma)=\aleph_2$, $(t,\Gamma,\name{\bar{\varphi}})\in{\bf G}(
\bA^q,\name{\cF}^q)$ and the type $\name{\gp}$ is (forced to be)
$\Gamma[\name{\bar{\varphi}}]$--big, \underline{then} there is a condition
$r\in\App$ as in (1) and such that $\name{\Gamma}_\gamma^r=\Gamma[
\name{\bar{\varphi}}]$.    
\end{enumerate}
\end{lemma}

\begin{proof}
1)\quad  Extend $\name{\cF}^q$ to $\name{\cF}^r$ so that $\name{x}_\gamma/
\name{\cF}^r$ realizes the required type.  
\medskip

\noindent 2) Note that every $\Gamma[\name{\bar{\varphi}}]$--big type can be 
extended to a complete $\Gamma[\name{\bar{\varphi}}]$--big one by \ref{1.1A}. 
\end{proof}

\begin{lemma}
\label{1.9Lemma}
\begin{enumerate}
\item Suppose $q_0,q_1,q_2\in\App$, $q_0=q_2\restriction\beta$, $q_0\leq
q_1$, $\bA^{q_1}\subseteq\beta$. Suppose further that $\bA^{q_2}\setminus
\bA^{q_0}= \{\beta\}$ and $\cf(\beta)=\aleph_2$. Assume that $\name{\gp}_1$
is a ${\mathbb C}_{\bA^{q_1}}$--name for a complete $\name{\Gamma}_\beta^{
q_2}$--big type over $(\mathop{{\prod^{\bA^{q_1}}}}\limits_{n<\omega}
\name{M}^n_{\bA^{q_1}}/\name{\cF}^{q_1})$ such that $\name{\gp}_1$ contains
the type $\name{\gp}_0$ induced by $\name{x}_\beta$ according to $q_2$. 
\underline{Then} there is $q_3\geq q_1,q_2$ with $\bA^{q_3}=\bA^{q_1}\cup
\{\beta\}$, such that $\name{x}_\beta$ induces $\name{\gp}_1$ on
$(\mathop{{\prod^{\bA^{q_1}}}}\limits_{n<\omega}\name{M}^n_{\bA^{q_1}}/ 
\name{\cF}^{q_1})$ (according to $q_3$).
\item Assume $q_0,q_1,q_2\in\App$, $q_0=q_2\restriction\beta$, $q_0\leq q_1$
and $\bA^{q_1}\subseteq\beta$. If $\bA^{q_2}\setminus\bA^{q_0}=\{\beta\}$
and $\cf(\beta)<\aleph_2$, \underline{then} there is $q_3\in\App$, $q_3\geq
q_1,q_2$ such that $\bA^{q_3}=\bA^{q_1}\cup\bA^{q_2}$. 
\item Assume that $\delta_1$, $\delta_2<\aleph_2$, and $\langle\beta_j:j<
\delta_2\rangle$ is a non-decreasing sequence of ordinals below $\aleph_3$.  
Let $\langle p_i:i<\delta_1\rangle$ be an increasing sequence from
$\App$. Suppose that $q_j\in \App\restriction\beta_j$ (for $j<\delta_2$) are
such that:
\[p_i\restriction\beta_j\leq q_j\mbox{ for }i<\delta_1,\ j<\delta_2,\qquad
q_j\leq_{{\rm end}}q_{j'}\mbox{ for }j<j'<\delta_2.\]
\underline{Then} there is an $r\in\App$ with $p_i\leq r$ and $q_j\leq_{{\rm
end}} r$ for all $i<\delta_1$ and $j<\delta_2$. 
\item If $\bar{p}=\langle p_i:i<\delta_1\rangle$ an increasing sequence in
$\App$, $\delta_1<\aleph_2$, \underline{then} $\bar{p}$ has an upper bound
in $\App$.  
\end{enumerate}
\end{lemma}

\begin{proof}
1)\quad Let $\bA_i=\bA^{q_i}$ and let $\name{\cF}_i=\name{\cF}^{q_i}$ for
$i<3$, and $\bA_3=\bA_1\cup\bA_2=\bA_1\cup\{\beta\}$. The only possibly
not clear part is to show that, in $\bV^{{\mathbb C}_{\bA_3}}$, there is an 
ultrafilter extending $\name{\cF}_1\cup\name{\cF}_2$ which contains
$\name{\cF}'$, the family of all the sets  
\[\{n<\omega:\name{M}^n_{\bA_3}\models\name{\varphi}[\name{x}_\beta(n),
\name{\bar{a}}(n)]\}\]
for $\name{\varphi}(x,\bar{y})\in\name{\gp}_1$, $\ell g(\bar{y})=m$, and a
${\mathbb C}_{\bA_1}$--name $\name{\bar{a}}$ of an $m$--tuple from
$\name{\omega}^*_{\bA_1}$ (and in our notation above $\name{\bar{a}}(n)$ is
a ${\mathbb C}_{\bA_1}$--name for an $m$--tuple of elements of $\omega$). As
$\name{\cF}_1,\name{\cF}_2,\name{\cF}'$ are (forced, i.e.,
$\forces_{{\mathbb C}_{\bA_3}}$) to be closed under intersections (of two, 
and hence of finitely many), clearly if this fails, then there are
$m<\omega$, a condition $p\in{\mathbb C}_{\bA_3}$, a ${\mathbb
C}_{\bA_1}$--name $\name{\bf a}$ of a member of $\name{\cF}_1$, a ${\mathbb
C}_{\bA_2}$--name $\name{\bf b}$ of a member of $\name{\cF}_2$, a (name for
a) $\tau_{\bA_1}$--formula $\name{\varphi}$ and a ${\mathbb
C}_{\bA_1}$--name for an $m$--tuple $\name{\bar{a}}$ from
$\name{\omega}^*_{\bA_1}$ such that 
\[p\restriction\bA_1\forces_{{\mathbb C}_{\bA_1}}\mbox{`` }\name{\varphi}(
x,\name{\bar{a}})\in\name{\gp}_1\mbox{ ''\quad and\quad }p\forces_{{\mathbb
C}_{\bA_3}}\mbox{`` }\name{\bf a}\cap\name{\bf b}\cap\name{\bf c}=\emptyset
\mbox{ '',}\]
where 
\[\name{\bf c}=\{n:\name{M}^n_{\bA_3}\models\name{\varphi}[
\name{x}_\beta(n), \name{\bar{a}}(n)]\}.\]
We may easily eliminate parameters, so we may assume that we have
$\name{\varphi}[\name{x}_\beta(n)]$ only (remember the definition of
$\tau_{\bA_1}$). Let $p_i=p\restriction \bA_i$ for $i=0,1,2$, and let
$H^0\subseteq{\mathbb C}_{\bA_0}$ be generic over $\bV$ such that $p_0\in
H^0$. For $n<\omega$ let $\name{A}^*_n$ be a ${\mathbb C}_{\bA_0}$--name
such that   
\[\begin{array}{ll}
\name{A}^*_n[H^0]=\{y\in\name{M}_{\bA_2}^n:&\mbox{there is } p_2'\in{\mathbb 
C}_{\bA_2}\mbox{ such that}\\ 
&p_2\leq p'_2,\ p'_2\restriction \bA_0
\in H^0\mbox{ and}\\
&p'_2\forces\mbox{`` }\name{x}_\beta(n)=y\mbox{ and }n\in\name{\bf b}
\mbox{ ''}\} 
  \end{array}\]
(recall $y\in \name{M}^n_{\bA_2}$ means $y\in\omega$). Let $\name{A}^*=
\mathop{{\prod}^{\bA_0}}\limits_{n<\omega}\name{A}^*_n/\name{\cF}_0$. So
$\name{A}^*[H^0]$ is (the interpretation of) an unary predicate from
$\tau_{\bA_0}$; in fact $Q_{\langle\name{A}^*_n:n<\omega\rangle}$ is such a
predicate. Thus, in $\bV[H_0]$, either $\name{A}^*(x)\in\name{\gp}_0$ or
$\neg\name{A}^*(x)\in\name{\gp}_0$. The latter is impossible by the choice
of $\name{\gp}_0$, so necessarily $\name{A}^*(x)\in\name{\gp}_0$. As also
$\forces_{{\mathbb C}_{\bA_1}}\mbox{`` }\name{\bar{\varphi}}(y)\in
\name{\gp}_1\mbox{ ''}$, clearly if $H^1\subseteq {\mathbb C}_{\bA_1}$ is
generic over $\bV$ and $H^0\cup\{p_1\}\subseteq H^1$, then in $\bV[H^1]$ we
have      
\[\{n\in\omega: \name{M}^n_{\bA_1}\models(\exists y)\big(\name{A}^*(y)\ \&\
\name{\varphi}(y)\big)\}\in\name{\cF}_1[H^1]\]
(remember $\name{\gp}_1$ is a type over $\mathop{{\prod^{ \bA^{q_1}}}}
\limits_{n<\omega}\name{M}^n_{\bA^{q_1}}/\name{\cF}_1$ extending
$\name{\gp}_0$). Consequently, we may find a condition $p'_1\in H_1\subseteq
{\mathbb C}_{\bA_1}$ stronger than $p_1$, an integer $n<\omega$, and an
element $y\in\name{M}_{\bA_1}^n$ (so $y\in \omega$) such that   
\[p'_1\restriction\bA_0\in H^0,\quad \mbox{ and }\quad p'_1\forces_{{\mathbb
C}_{\bA_1}}\mbox{`` }\name{M}^n_{\bA_1}\models\big(\name{A}^*(y)\ \&\
\name{\varphi}(y)\big)\mbox{ and }n\in \name{\bf a}\mbox{ ''}.\] 
As $\name{A}^*_n$ is a ${\mathbb C}_{\bA_0}$--name, we really have $y\in
\name{A}_n^*[H^0]$, and hence (by its definition) for some $p'_2\in{\mathbb
C}_{\bA_2}$ we have 
\[p_2\leq p'_2,\quad p'_2\restriction\bA_0\in H^0,\quad\mbox{ and }\quad
p'_2 \forces\mbox{`` }y=\name{x}_\beta(n)\mbox{ and }n\in\name{\bf b}
\mbox{ ''}.\]
Now for our $n$ we can force $n\in\name{\bf a}\cap\name{\bf b}\cap\name{\bf
c}$ by amalgamating the corresponding conditions $p'_1,p'_2$, getting a
contradiction. As said above this finishes the proof of the existence of
$q_3$. 
\medskip

\noindent 2)\quad The proof is essentially contained in the previous one
(use the very trivial bigness notion: $\varphi(x,\bar{a})$ is big in $M$ if
and only if $M\models(\exists x)\varphi(x,\bar{a})$, so we may use a
$\name{\gp}_1$). See also the end of the proof of (3). 
\medskip

\noindent 3)\quad We will prove by induction on $\gamma\in\aleph_3$ that if
all $\beta_j\leq \gamma$ and all $p_i$ belong to $\App\restriction\gamma$,
then the assertion in (3) holds for some $r\in\App\restriction\gamma$. 

We may assume that $\delta_1>0$ (otherwise apply \ref{1.7Lemma}) and
$\delta_2>0$ (otherwise let $\delta_2'=1$, $\beta_0=0$, $q_0'\in\App
\restriction 0$ be above $p_i\restriction 0$ for $i<\delta_1$; so it just
means $\name{\cF}^{q_0'}$ is the ultrafilter $\name{\cF}^{p_i\restriction
0}$ for $i<\delta_1$; now if $\gamma=0$, then $r=q_0'$ is as required and
otherwise we have reduced the case $\delta_2=0$ to the case $\delta_2=1$). 

We may assume that $\beta_j=\sup\{\alpha+1:\alpha\in \bA^{q_j}\}$ (for
$j<\delta_2$), and also that the sequence $\langle\beta_j:j<\delta_2\rangle$
is strictly increasing. Let $\beta=\sup\limits_{j<\delta_2}\beta_j$ and let 
$q=(\bigcup\limits_{j<\delta_2}\bA^{q_j},\bigcup\limits_{j<\delta_2}
\name{\cF}^{q_j},\bigcup\limits_{j<\delta_2}\name{\bar{\Gamma}}^{q_j})$. 
\smallskip

We first assume $\cf(\gamma)\neq\aleph_0$. 
\smallskip

\noindent If $\gamma=\beta$, then $q\in\App$ and we may take $r=q$. So let
us assume $\beta<\gamma$. If $\delta_2$ is a successor ordinal, or a limit
ordinal of uncountable cofinality, then we let $q^*=q$ (clearly $q^*\in\App
\restriction\beta$). If $\cf(\delta_2)=\aleph_0$, then we may first apply
the inductive hypothesis to $\langle p_i\restriction\beta:i<\delta_1\rangle$
(and $\langle\beta_j,q_j:j<\delta_2\rangle$) to get a condition $q^*\in\App
\restriction\beta$ which is stronger than all $p_i\restriction\beta$ and
which end-extends all $q_j$. So in all these cases, we have a condition $q^*
\in\App\restriction\beta$ (end extending all $q_j$ for $j<\delta_2$)
stronger than all $p_i\restriction\beta$ for $i<\delta_1$ (and we are
looking for its end-extension which is a bound to all $p_i\restriction
\beta$).
\smallskip

\noindent{\sc The case $\gamma=\gamma_0+1$, a successor}\\
In this case our inductive hypotheses applies to the $p_i\restriction
\gamma_0,q^*$, and $\gamma_0$, yielding $r_0$ in $\App\restriction\gamma_0$
with $p_i\restriction\gamma_0\leq r_0$ and $q^*\leq_{{\rm end}}r_0$. What
remains to be done is an amalgamation of $r_0$ with all of the $p_i$, where
$\bA^{p_i}\subseteq\bA^{r_0}\cup\{\gamma_0\}$, and where one may as well
suppose that $\gamma_0$ is in $\bA^{p_i}$ for all $i$. This is a slight
variation on (1) or (2). For instance, suppose $\cf(\gamma_0)=\aleph_2$. We 
let  
\begin{itemize}
\item $\bA_2=\bigcup\limits_{i<\delta_1}\bA^{p_i}$, $\bA_0=\bA_2\setminus
\{\gamma_0\}$, $\bA_1=\bA^{r_0}$, $\bA_3=\bA_2\cup\bA_1$.
\item $\name{\cF}_1=\name{\cF}^{r_0}$, $\name{\cF}_2=\bigcup\limits_{i<
\delta_1}\name{\cF}^{p_i}$ (the latter might be only a ${\mathbb
C}_{\bA_2}$--name of a filter).
\item For $i<\delta_1$ let $\name{\gp}^i$ be the ${\mathbb C}_{\bA^{p_i}\cap 
\gamma_0}$--name for the ($\name{\Gamma}^{p_i}_{\gamma_0}$--big) type 
induced by $\name{x}_{\gamma_0}$ over the model $\mathop{{\prod}^{\bA^{p_i}
\cap\gamma_0}}\limits_{n<\omega}\name{M}^n_{\bA^{p_i}\cap\gamma_0}/
\name{\cF}^{p_i\restriction\gamma_0}$. Then let $\name{\gp}_0=
\bigcup\limits_{i<\delta_1}\name{\gp}^i$, and note that it is a ${\mathbb
C}_{\bA_0}$--name for a $\name{\Gamma}^{p_i}_{\gamma_0}$--big type over the
model $\mathop{{\prod}^{\bA_0}}\limits_{n<\omega}\name{M}^n_{\bA_0}/
\name{\cF}_0$. 
\item Let $\name{\gp}_1$ be (a ${\mathbb C}_{\bA_1}$--name for) a complete 
$\name{\Gamma}_{\gamma_0}^{p_i}$--big type over $\mathop{{\prod}^{\bA_1}}
\limits_{n<\omega}\name{M}^n_{\bA_1}/\name{\cF}_0$ extending
$\name{\gp_0}$. (Exists by \ref{1.1A}.)
\end{itemize}
Now, in $\bV^{{\mathbb C}_{\bA_3}}$, we want to extend $\name{\cF}_1\cup
\name{\cF}_2$ to an ultrafilter $\cF'$ containing the sets of the form
$\{n<\omega:\name{M}^n_{\bA_3}\models\name{\varphi}[\name{x}_{\gamma_0}(n)]
\}$ for all $\name{\varphi}(x)\in\name{\gp}_1$. If this fails, then as 
\[\forces_{{\mathbb C}_{\bA_1}}\mbox{`` }\langle \name{\cF}^{p_i}:i<\delta_1
\rangle\mbox{ is increasing ''}\]
we find a condition $p\in{\mathbb C}_{\bA_3}$, a ${\mathbb C}_{\bA_1}$--name
$\name{\bf a}$ of a member of $\name{\cF}_1$, and $i<\delta_1$, and a
${\mathbb C}_{\bA_2}$--name $\name{\bf b}$ for a member of $\name{\cF}_i$,
and $\name{\varphi}$ such that  
\[p\restriction\bA_1\forces\name{\varphi}(x)\in\name{\gp}^i\subseteq
\name{\gp}_1\qquad\mbox{and}\qquad p\forces_{{\mathbb C}_{\bA_3}}\name{\bf
a}\cap\name{\bf b}\cap\{n:M^n_{\bA_3}\models\name{\varphi}[x_\beta(n)]\}=
\emptyset.\]
Next we continue exactly as in the proof of (1).
\smallskip

\noindent{\sc The case $\gamma$ is a limit ordinal of cofinality
$\aleph_2$}\\  
Since $\delta_1<\aleph_2$ there is some $\gamma_0<\gamma$ such that all
$p_i$ lie in $\App\restriction\gamma_0$ and $\beta<\gamma_0$, and the
induction hypothesis then yields the claim. 
\smallskip

\noindent{\sc The case $\gamma$ is a limit ordinal of cofinality
$\aleph_1$}\\ 
Choose a strictly increasing and continuous sequence $\langle\gamma_j: j<
\aleph_1\rangle$ with supremum $\gamma$, starting with $\gamma_0=\beta$.
By induction on $j$ choose $r_j\in\App\restriction\gamma_j$ (for $j<
\aleph_1$) such that:
\begin{itemize}
\item $r_0=q^*$;
\item $r_j\leq_{{\rm end}}r_{j'}$ for $j<j'<\aleph_1$;
\item $p_i\restriction\gamma_j\leq r_j$ for $i<\delta_1$ and $j<\aleph_1$.
\end{itemize}
[Thus, at a successor stage $j+1$, the inductive hypothesis is applied to
$p_i\restriction\gamma_{j+1},r_j,\gamma_j$, and $\gamma_{j+1}$. At a limit
stage $j$, we apply the inductive hypothesis to $p_i\restriction\gamma_j$
for $i<\delta_1$, $r_{j'}$ for $j'<j$, $\gamma_{j'}$ for $j'<j$, and
$\gamma_j$.] Finally, we let $r=(\bigcup\limits_{j<\aleph_1}\bA^{r_j},
\bigcup\limits_{j<\aleph_1}\name{\cF}^{r_j},\bigcup\limits_{j<\aleph_1}
\name{\bar{\Gamma}}^{r_j})$. Clearly $r\in\App$ is as required.   
\medskip

Now we are going to consider the remaining case:
\smallskip

\noindent{\sc The case $\gamma$ is a limit ordinal of cofinality
$\aleph_0$}\\ 
If $\beta<\gamma$ (where $\beta$ is as defined at the beginning of the
proof), then we first pick a strictly increasing sequence $\langle\gamma_j: 
j<\aleph_0\rangle$ of ordinals such that $\beta<\gamma_0$ and
$\sup\limits_{j<\aleph_0}\gamma_j=\gamma$. Then we apply repeatedly the
inductive hypothesis to build a sequence $\langle q_j':j<\aleph_0\rangle$
such that $q_j'\in\App\restriction\gamma_j$, $q_{j_0}'\leq_{\rm end}
q_{j_1}'$ for $j_0<j_1$, $q_j\leq_{\rm end} q_0'$ (for all $j<\aleph_0$),
and $p_i\restriction\gamma_j\leq q_j'$ (for all $i<\delta_1$,
$j<\aleph_0$). Thus we have reduced this sub-case to the only one remaining:
$\beta=\gamma$. Now if for some $j<\delta_2$ we have $\beta_j=\gamma$, then
$r=q_j$ is as required, so without loss of generality $(\forall j<\delta_2)
(\beta_j<\gamma)$. Then necessarily $\cf(\delta_2)=\aleph_0$ and we may
equally well assume that $\delta_2=\aleph_0$. 

We take $q$ as defined earlier (so it is the ``union'' of all $q_j$), but
it does not have to be a condition in $\App$: the filter $\bigcup\limits_{j< 
\aleph_0}\name{\cF}^{q_j}$ does not have to be an ultrafilter, and we need
to extend it to one that contains also $\bigcup\limits_{i<\delta_1}
\name{\cF}^{p_i}$. Note that $\bA^*\stackrel{\rm def}{=}\bigcup\limits_{i<
\delta_1}\bA^{p_i}\subseteq\bigcup\limits_{j<\aleph_0}\bA^{q_j}
\stackrel{\rm def}{=}\bA^+$, but there might be ${\mathbb C}_{\bA^*}$--names
for elements of $\bigcup\limits_{i<\delta_1}\name{\cF}^{p_i}$ that are not
${\mathbb C}_{\bA^{q_j}}$--names for any $j<\aleph_0$, so it could happen
that one name like that is forced to be disjoint from some element of
$\name{\cF}^{q_j}$. So assume toward contradiction, that there are a
condition $p\in {\mathbb C}_{\bA^+}$, ordinals $i<\delta_1$ and
$j<\aleph_0$, a ${\mathbb C}_{\bA^{p_i}}$--name $\name{\bf a}$, and a
${\mathbb C}_{\bA^{q_j}}$--name $\name{\bf b}$ such that
\[p\forces_{{\mathbb C}_{\bA^+}}\mbox{`` }\name{\bf a}\in\name{\cF}^{p_i}\ 
\&\ \name{\bf b}\in\name{\cF}^{q_j}\ \&\ \name{\bf a}\cap\name{\bf b}=
\emptyset\mbox{ ''.}\]
Increasing $j$ if necessary, we may also assume that $p\in {\mathbb
C}_{\bA^{q_j}}$ so $\Dom(p)\subseteq\beta_j\times\omega$. Let $H^0\subseteq
{\mathbb C}_{\bA^{p_i}\cap\beta_j}$ be generic over $\bV$ such that
$p\restriction\bA^{p_i}\in H^0$, and let    
\[\begin{array}{ll}
{\bf c}=\{n\in\omega: &\mbox{there is a condition }p'\in {\mathbb C}_{
\bA^{p_i}}\mbox{ stronger than }p\restriction\bA^{p_i}\mbox{ and}\\
&\mbox{such that }p'\restriction (\bA^{p_i}\cap\beta_j)\in H^0\ \mbox{and }\
p'\forces_{{\mathbb C}_{\bA^{p_i}}}\mbox{`` } n\in \name{\bf a}\mbox{ ''}\}.
  \end{array}\]  
Clearly, ${\bf c}\in\bV[H^0]$ is a set from $\big(\name{\cF}^{p_i}
\restriction (\bA^{p_i}\cap \beta_j)\big)[H^0]$. Since $p_i\restriction
\beta_j\leq q_j$, we find a condition $p''\in {\mathbb C}_{\bA^{q_j}}$ and
$n\in{\bf c}$ such that  
\[p\leq p''\quad \&\quad p''\restriction (\bA^{p_i}\cap\beta_j)\in H^0\quad
\&\quad p''\forces_{{\mathbb C}_{\bA^{q_j}}}\mbox{`` }n\in \name{\bf
b}\mbox{ ''}.\] 
For this $n$ we take $p'\in {\mathbb C}_{\bA^{p_i}}$ witnessing that $n\in
{\bf c}$ and next we let $p^*=p'\cup p''$. Clearly $p^*\forces n\in\name{\bf 
a}\cap \name{\bf b}$, a contradiction.
\medskip

\noindent 4)\quad Follows, i.e., it is the case $\delta_2=0$ of part (3).
\end{proof}

\begin{lemma}
\label{cc}
Assume $\bV\models{\rm GCH}$. The forcing notion $\App$ satisfies the
$\aleph_3$--chain condition, it is $\aleph_2$--complete, $\|\App\|=\aleph_3$
and $\|\App\restriction\gamma\|\leq \aleph_2$ for every $\gamma\in\aleph_3$. 
Consequently, the forcing with $\App$ does not collapse cardinals nor
changes cofinalities, and $\forces_{\App}{\rm GCH}$. 
\end{lemma}

\begin{proof}
The only perhaps unclear part is the chain condition. Suppose we have an
antichain $\{q_\alpha:\alpha\in\aleph_3\ \&\ \cf(\alpha)=\aleph_2\}\subseteq
\App$ (the index $\alpha$ is taken to vary over ordinals of cofinality
$\aleph_2$ just for convenience). An important point is that $\bG$ can
``offer'' at most $\aleph_2$ candidates for the bigness notion at
$\delta<\aleph_3$, $\cf(\delta)=\aleph_2$, hence for each $\gamma\in
\aleph_3$ the restricted forcing $\App\restriction\gamma$ has cardinality
$\leq\aleph_2$. Applying Fodor's lemma twice, we find a stationary set
$S\subseteq\{\alpha\in\aleph_3:\cf(\alpha)=\aleph_2\}$ and a condition
$q^*\in\App$ such that $(\forall\alpha\in S)(q_\alpha\restriction\alpha=
q^*)$. Pick $\alpha_1,\alpha_2\in S$ such that $\sup(\bA^{q_{\alpha_1}})<
\alpha_2$; it follows from Lemma \ref{1.9Lemma} that the conditions
$q_{\alpha_1},q_{\alpha_2}$ are compatible, a contradiction.  
\end{proof}

\begin{proposition}
\begin{enumerate}
\item For each $p\in\App$ and $\alpha\in\aleph_3$, there is a condition
$q\in\App$ stronger than $p$ and such that $\alpha\in\bA^{q}$.
\item $\name{\cF}\stackrel{\rm def}{=}\bigcup\{\name{\cF}^r: r\in
\name{G}_\App\}$ is a $\coh$--name for an ultrafilter on $\omega$. Also, for
each $r\in \name{G}_\App$ we have: $\name{\cF}\cap{\mathcal
P}(\omega)^{(\bV^*)^{{\mathbb C}_{\bA^r}}}=\name{\cF}^r$.   
\end{enumerate}
\end{proposition}

\begin{proof}
Should be clear (for (1) use \ref{1.8Lemma} + \ref{1.9Lemma}(1); then (2)
follows).  
\end{proof}

\begin{definition}
\label{2.8}
\begin{enumerate}
\item Suppose $G_\App\subseteq\App$ is generic over $\bV$, $\bV^*=\bV[
G_\App]$. We let $\name{\cF}^\delta$ be the ${\mathbb C}_\delta$--name for
the restriction $\name{\cF}\restriction\delta$ of the ultrafilter
$\name{\cF}$ to the sets from the universe $(\bV^*)^{{\mathbb C}_\delta}$. 
\item We define an $\App$--name $\name{\Gamma}_\delta$ of a ${\mathbb
C}_\delta$--name as $\Gamma^p_\delta$ for every $p\in\name{G}_\App$ such
that $\delta\in\bA^p$. (So it is an $\App*{\mathbb C}_\delta$--name.)
\end{enumerate}
\end{definition}

\begin{lemma}
\label{1.12Lemma}
\begin{enumerate}
\item Suppose that $G_\App\subseteq\App$ is generic over $\bV$, $\bV^*=\bV[ 
G_\App]$, and $\delta<\aleph_3$, $\cf(\delta)=\aleph_2$, and $H^\delta
\subseteq {\mathbb C}_\delta$ is generic over $\bV^*$. \underline{Then}, in
$\bV[G_\App\cap(\App\restriction\delta)][H^\delta]$, we have\footnote{Note:
$\name{M}_\delta^n$ is $\name{M}^n_\bA$ for $\bA=\delta$}:    
\[\prod\limits_{n<\omega}\name{M}^n_\delta/\name{\cF}^\delta[H^\delta]\
\mbox{ is $\aleph_2$--compact.}\]
\item Also if $H\subseteq {\mathbb C}_{\aleph_3}$ is generic over $\bV^*$,
$H\supseteq H^\delta$, \underline{then} in $\bV^*[H]$:
\begin{enumerate}
\item[(a)] $\prod\limits_{n<\omega} M^n_{\aleph_3}/\name{\cF}[H]$ is
$\aleph_2$--compact,
\item[(b)] $\name{x}_\delta[H]/\name{\cF}[H]\in \prod\limits_{n<\omega}
M^n_{\aleph_3}/\name{\cF}[H]$ realizes a $\name{\Gamma}_\delta[G]
[H^\delta]$--big type over $\mathop{\prod^\delta}\limits_{n<\omega}
M^n_\delta/\name{\cF}^\delta[H^\delta]$. 
\end{enumerate}
\end{enumerate}
\end{lemma}

\begin{proof}
By \ref{1.8Lemma}(1). We can use some $\name{x}_\beta$ with $\beta$ of
cofinality less than $\aleph_2$ to realize each type.  
\end{proof}

\section{Definability}

\begin{hypothesis}
\label{hypo}
In this section we assume that $\bG$ is an $(\aleph_3,\aleph_2)$--bigness
guide, $\App=\App_\bG$, $G\subseteq\App$ is a generic filter over $\bV$, and
$\bV^*=\bV[G]$. For an ordinal $\delta<\aleph_3$, we let $G_\delta=G\cap
(\App\restriction\delta)$. Also, $\name{H}$, $\name{H}^\delta$ are the
canonical $\coh$-- and ${\mathbb C}_\delta$--names of the generic subsets of
${\mathbb C}_{\aleph_3}$ and ${\mathbb C}_\delta$, respectively. We work
mostly in $\bV^*$.  
\end{hypothesis}
\noindent [Note that, by Lemma \ref{cc}, $\bV^*\models{\rm GCH}$.]

\begin{definition}
\label{3.1B}
\begin{enumerate}
\item We say that $\bm$ is {\em an $(\aleph_3,\aleph_2)$--isomorphism
candidate} (or just {\em an isomorphism candidate}, in $\bV$ or in $\bV^*$,
see below; letting $\bm^-=\langle \name{t},\name{\bar{\varphi}},\name{\psi},
\name{\Delta},\langle \name{N}^\ell_n:n<\omega,\ \ell=1,2\rangle\rangle$,
note that as $\App$ is $\aleph_2$--complete, this forcing does not add new
$\bm^-$, i.e., $\bV$ and $\bV^*$ have the same set of $\bm^-$, though we
have an $\App$--name $\name{\bm}$ of such object) if:
\begin{enumerate}
\item[(i)] $\bm$ consists of $\bA^*=\bA^*[\bm]\subseteq [\aleph_3]^{<
\aleph_2}$, $P^*=P^*[\bm]$, $\name{N}^\ell_n=\name{N}^\ell_n[\bm]$ (for $n<
\omega$, $\ell\in\{1,2\}$), $\name{F}=\name{F}[\bm]$, $\name{\Gamma}=
\name{\Gamma}[\bm]$ and $(\name{t},\name{\bar{\varphi}},\name{\psi},
\name{\Delta})=(\name{t}[\bm],\name{\bar{\varphi}}[\bm],\name{\psi}[\bm],
\name{\Delta}[\bm])$;
\item[(ii)] $\name{t},\name{\bar{\varphi}},\name{\psi},\name{\Delta}$ are
${\mathbb C}_{\bA^*}$--names as in \ref{2.1A}(3) and $\name{\Gamma}=
\name{\Gamma}_{(\name{\bar{\varphi}},\name{t},\name{\psi})}$ is a bigness
notion as there, $\tau(\name{t})$ countable for simplicity; 
\item[(iii)] $\name{N}_n^\ell$, for $n<\omega$ and $\ell\in\{1,2\}$, are
${\mathbb C}_{\bA^*}$--names for countable models of a (countable) theory
$t^\ell_n$, and the universes $|\name{N}_n^\ell|$ are subsets of $\omega$. 

\noindent Also it is forced (i.e., $\forces_{\coh}$) that $\name{t}=
{\rm Th}\Big(\prod\limits_{n<\omega}\name{N}^\ell_n/\name{\cF}\Big)={\rm Th}
\Big(\prod\limits_{n<\omega} N^\ell_n/\name{\cF}\Big)$, so the
$\prod\limits_{n<\omega}$ is $\mathop{{\prod}^{\aleph_3}}\limits_{n<
\omega}$.

\item[(iv)] We have predicates $Q_R^\ell\in\tau_{\bA^*}$ (for $R\in\tau(t)$)
such that $\name{\bar{\varphi}}^\ell=\langle Q_R^\ell:R\in\tau(t)\rangle$ is
the interpretation of $\name{t}$ in $\prod\limits_{n<\omega}
\name{M}_{\bA^*}^n/\name{\cF}$ giving $\prod\limits_{n<\omega}
\name{N}^\ell_n/\name{\cF}$. (Remember \ref{generalia}(4),
\ref{verygoodb}(1); so by the choice of $\tau_\bA$ actually
$\name{\bar{\varphi}}^*=\bar{\varphi}^*$.) 
\item[(v)] $\name{F}$ is a $\coh$--name (more accurately an $\App$--name of
such name, but we sometimes write $\name{\cF}$ instead of $\name{\cF}[G]$
as when $G$ is constant) and $p^*\in\coh$ is a condition such that:     
\[\begin{array}{ll}
\forces_{\App}\ \ p^*\forces_{\coh}&\mbox{`` }\name{F}\mbox{ is a map from 
}\prod\limits_{n<\omega}\name{N}^1_n\mbox{ into }\prod\limits_{n<\omega}
\name{N}_n^2\\
&\quad\mbox{which represents a $\name{\Delta}$--embedding modulo }\name{\cF}
\mbox{ ''.} 
  \end{array}\]
\end{enumerate}
\item $\bm$ is {\em $(P,\xtheta)$--separative} if $P,\xtheta$ are ${\mathbb
C}_{\bA^*}$--names and there is a witness $\name{X}\subseteq
\omega^*_{\bA^*}$ for $\name{\Gamma}^{\bm}$ in the intended model, i.e.,
this is forced, $\forces_{{\mathbb C}_{\bA^*}}$. 
\end{enumerate}
[If $\bm$ is clear from the context we may omit it.]
\end{definition}

\begin{observation}
\label{3.1C}
Assume, in $\bV$, that $\bm$ is an $(\aleph_3,\aleph_2)$--isomorphism
candidate, $\name{\Gamma}=\name{\Gamma}[\bm]=\Gamma_{(\name{t},
\name{\bar{\varphi}},\name{\psi})}$.
\underline{Then} there is a stationary set of ordinals $\delta<\aleph_3$
such that: 
\begin{enumerate}
\item[(a)$_\delta$] $\bA^*\subseteq\delta$, $\cf(\delta)=\aleph_2$, and
$p^*\in\coh\restriction\delta$, and for some $q\in G$ we have that
$\name{\Gamma}^q_\delta$ is $\Gamma[\name{\bar{\varphi}}]$ (for $(t,\Gamma,
\name{\bar{\varphi}})$ from \ref{2.1A}(3)),
\item[(b)$_\delta$]  for every $\coh\restriction\delta$--name $\name{x}$ for
an element of $\prod\limits_{n<\omega}\name{N}^1_n$, $\name{F}(\name{x})$ is
a $\coh\restriction\delta$--name,  

\noindent [recall $\App$ satisfies the $\aleph_3$--cc]
\item[(c)$_\delta$] similarly for $\name{F}^{-1}$ and for `` $y\in\Rang(
\name{F})$ '',
\item[(d)$_\delta$] $\forces_\coh\{n<\omega:\name{x}_\delta(n)\in
\name{N}^1_n\}\in\name{\cF}$ (so $\name{x}_\delta/\name{\cF}\in
\prod\limits_{n<\omega}\name{N}^1_n/\name{\cF}$).
\end{enumerate}
For such $\delta$, we let $\name{y}^*=\name{y}^*_\delta=\name{y}^*_{\delta,
\name{F}}=y^*_{\delta,\bm}$ be $\name{F}(\name{x}_\delta)\in\prod\limits_{n<
\omega}\name{N}_n^2$. 
\end{observation}

\begin{mic}
\label{3.1D}
Assume that $\bm$ is an $(\aleph_3,\aleph_2)$--isomorphism candidate as in
\ref{3.1C}, and $\delta<\aleph_3$ is as there. \underline{Then} there are
$q_\delta,y$ such that    
\begin{enumerate}
\item[(a)] $q_\delta\in\App$, moreover $q_\delta\in G$,
\item[(b)] $q_\delta\forces_\App\mbox{`` }\name{F}(x_\delta)=\name{y}^*
\mbox{ ''}$, where $\name{y}^*$ is a ${\mathbb C}_{\bA^{q_\delta}}$--name of
a member of ${}^\omega\omega$, 
\item[(c)] $\bA^*\subseteq \bA^{q_\delta}$, $\bA_\delta\stackrel{\rm def}{=}
\bA^{q_\delta}\cap\delta$,
\item[(d)] in $\bV[G_\delta][\name{H}^\delta]$ we have:
\begin{enumerate}
\item[(i)] $\cF_\delta=\name{\cF}_\delta[G_\delta][\name{H}^\delta]$ is
a non-principal ultrafilter on $\omega$.
\item[(ii)] The model $M_\delta=\mathop{{\prod}^\delta}\limits_{n<\omega}
M^n_\delta/\cF_\delta$ with the vocabulary $\tau_\delta$ is
$\aleph_2$--compact. 
\item[(iii)] The vocabulary $\tau_{\bA_\delta}\subseteq\tau_\delta$ is of
cardinality $\leq\aleph_1$.
\item[(iv)] $M_{\bA_\delta}=\mathop{{\prod}^{\bA_\delta}}\limits_{n
<\omega} M^n_{\bA_\delta}/\cF^{q_\delta\restriction\delta}[\name{H}^\delta]
\prec M_\delta\restriction\tau_{\bA_\delta}$.  
\item[(v)] $\name{F}_\delta=(\name{F}\restriction\delta)[\name{H}^\delta]=
\big((\name{F}\restriction\delta)[G\cap(\App\restriction\beta)]\big)[
\name{H}^\delta]$ is a $\name{\Delta}$--embedding from the model
$\mathop{{\prod}^\delta}\limits_{n<\omega} N^1_n/\cF_\delta$ into
$\mathop{{\prod}^\delta}\limits_{n<\omega}N^2_n/\cF_\delta$,  
\item[(vi)]    Let $\name{p}_\delta=\name{p}_\delta(x)$ be the (${\mathbb
C}_\delta$--name of the) 1--type in the vocabulary $\tau_{\bA_\delta}$ such
that  
\[\mbox{`` }p_\delta(x)\mbox{ is the type realized by $\name{x}_\delta$ over
$M_{\bA_\delta}$ in $\mathop{{\prod}^{\bA_\delta}}\limits_{n<\omega}
M^n_{\bA^{q_\delta}}/\name{\cF}^{q_\delta}$ ''}.\] 
[Clearly it is a ${\mathbb C}_\delta$--name, or an $\App*{\mathbb
C}_\delta$--name; see clause (d) of Definition \ref{app}(1).]\\
\underline{Then} $p_\delta$ is $\Gamma$--big.
\item[(vii)]   Let $N^\ell_\delta=\mathop{\prod^\delta}\limits_{n<\omega}
N^\ell_n/\name{\cF}_\delta$ (they are in $\bV^*[\name{H}^\delta]$, even in
$\bV[G \restriction\delta][\name{H}^\delta]$). We define
$R_{\delta,m}\subseteq (N^1_\delta)^m\times (N^2_\delta)^m$ for $m<\omega$
so that 
\begin{enumerate}
\item[$(\circledast)_1$] $R_{\delta,m}$ includes the graph of $F_\delta$,
i.e., if $\bar{a}$ is an $m$--tuple from $N^1_\delta$, then $(\bar{a},
F_\delta(\bar{a}))\in R_{\delta,m}$, 
\item[$(\circledast)_2$] the truth value of $(\bar{a},\bar{b})\in R_{\delta,
m}$ depends only on $L_{\omega,\omega}(\tau_{\bA_\delta})$--type
realized by $(\bar{a},\bar{b})$ over $M_{\bA_\delta}$ in $M_\delta$,
\item[$(\circledast)_3$] in fact $R_{\delta,m}$ is minimal such that
$(\circledast)_1+(\circledast)_2$ hold.
\end{enumerate}
\item[(viii)] The relations $R_{\delta,m}$ mentioned above satisfy:  
\begin{enumerate}
\item[$(\oplus)_1$]  if $\bar{a}_1,\bar{a}_2$ are finite sequences of the
same length $m$ of members of $N^1_\delta$, and $p_\delta\cup\{
\xtheta^{N_1}(x,\bar{a}_1),\neg\xtheta^{N_1}(x,\bar{a}_2)\}$ is a
$\Gamma$--big type over $M_\delta$, and $\xtheta,\neg\xtheta\in\name{\Delta}
[\bm]$, where $\xtheta^{N_1}$ is $\xtheta$ as interpreted in the
interpretation $\bar{\varphi}^1$,\\
\underline{then} $(\bar{a}_1,F_\delta(\bar{a}_2))\notin R_{\delta,m}$. 
\item[$(\oplus)_2$] Above, we may replace $\xtheta,\neg\xtheta$ by any pair
$\xtheta_0,\xtheta_1$ of contradictory formulas from $\name{\Delta}[\bm]$. 
\end{enumerate}
\item[(ix)] Note that also 
\[\hskip -.5cm
\begin{array}{ll}
(*)^{p^*}_{\name{y}^*,\delta}\quad  p^*\forces_{\coh}&\mbox{``the
$\name{\Delta}$--type which }\name{y}^*\mbox{ realizes over
}\name{N}^2_\delta=(\prod\limits_{n<\omega}\name{N}_n^2/\name{\cF})^{(
\bV^*)^{\coh\restriction\delta}}\\  
&\ \mbox{ in the model $\name{N}^2=(\prod\limits_{n<\omega}\name{N}_n^2/
\name{\cF})^{(\bV^*)^\coh}$ includes the image}\\
&\ \mbox{ under }\name{F}\mbox{ of the $\name{\Delta}$--type which}\\  
&\ \ \name{x}_\delta/\name{\cF}\mbox{ realizes over }\name{N}^1_\delta= 
(\prod\limits_{n<\omega}\name{N}^1_n/\name{\cF})^{(\bV^*)^{\coh\restriction 
\delta}}\\
&\ \mbox{ in the model $\name{N}^1=(\prod\limits_{n<\omega}
\name{N}_n^1/\name{\cF})^{(\bV^*)^\coh}$ ''.}
  \end{array}\]
\item[(x)] If $\name{\Delta}$ is closed under negation we can naturally note
that there are equivalence relations $E_\ell$ on $N^\ell_\delta$ satisfying
$(\circledast)_2$ induced by $(*)^{p^*}_{\name{y}^*,\delta}$. 
\end{enumerate}
\end{enumerate}
\end{mic}

Also notice that the clauses (b)$_\delta$, (c)$_\delta$ of \ref{3.1C} above
say that $\name{F}^\delta[G]$ is really a ${\mathbb C}_\delta$--name for a
function from $(\prod\limits_{n<\omega} \name{N}^1_n)^{(\bV^*)^{{\mathbb
C}_\delta}}$ into $(\prod\limits_{n<\omega} \name{N}^2_n)^{(\bV^*)^{{
\mathbb C}_\delta}}$ preserving $\name{\Delta}$--formulas; in the main case
it is ``onto''. 

\subsection*{The proof of the Main Isomorphism Theorem} 
The proof of \ref{3.1D} is broken into several steps and lemmas. Note that
we use the countability of $t$.

Take a condition $q_\delta\in G$ such that  
\begin{enumerate}
\item[(A)$^{q_\delta}$] $\bA^*\subseteq \bA^{q_\delta}$, $\name{x}_\delta,
\name{y}^*$ are ${\mathbb C}_{\bA^{q_\delta}}$--names (so $\delta\in
\bA^{q_\delta}$), and $p^*\in {\mathbb C}_{\bA^{q_\delta}}$, and  
\item[(B)$^{q_\delta}$] the condition  $q_\delta$ forces (in $\App$) that
clauses (b)$_\delta$, (c)$_\delta$ and (d)$_\delta$ from \ref{3.1C} hold
true (so in particular $q_\delta$ forces that $\name{x}_\delta/\name{\cF}\in 
\prod\limits_{n<\omega}\name{N}^1_n/\name{\cF}$, $\name{y}^*\in
\prod\limits_{n<\omega}\name{N}_n^2$ and $(*)^{p^*}_{\name{y}^*,\delta}$
holds), and 
\item[(C)$^{q_\delta}$] {\bf if} $\name{x}$ is a ${\mathbb C}_{
\bA^{q_\delta}}$--name for a member of $\mathop{{\prod}^{\bA^{q_\delta}}}
\limits_{n<\omega}\name{N}^1_n$ ($\mathop{{\prod}^{\bA^{q_\delta}}}
\limits_{n<\omega}\name{N}^2_n$, respectively), {\bf then}
$\name{F}(\name{x})$ ($\name{F}^{-1}(\name{x})$, respectively) is also a
${\mathbb C}_{\bA^{q_\delta}}$--name. 

\noindent [In fact, also $\name{X}\subseteq\mathop{{\prod}^{\bA^*}}\limits_{ 
n<\omega} M^n_{\bA^*}$ by \ref{3.1C}.]
\end{enumerate}
Before we continue with the proof of \ref{3.1D}, let us note the following. 

\begin{lemma}
\label{1.15Proposition}
Let $\delta<\aleph_3$, $q_\delta\in\App$ and $\name{y}^*,p^*$ be as above. 
Suppose that  
\[q_\delta\restriction \delta=q\leq q'\in G\cap(\App\restriction\delta).\] 
Let $\ytheta$ be a ${\mathbb C}_{\bA^*}$--name of a $\tau(t)$--formula. 
Assume further that $\name{x}'$, $\name{x}''$ and $\name{y}'$, $\name{y}''$
are ${\mathbb C}_{\bA^{q'}}$--names, and $p^*\leq p\in\coh$, and the
condition $p\restriction\bA^{q'}$ forces (in ${\mathbb C}_{\bA^{q'}}$) that   
\begin{enumerate}
\item[$(\alpha)$] $\name{x}',\name{x}''\in\prod\limits_{n<\omega}
\name{N}^1_n$, and $\name{y}',\name{y}''\in\prod\limits_{n<\omega} 
\name{N}_n^2$, and   
\item[$(\beta)$]  the types of $(\name{x}',\name{y}')$ and of $(\name{x}'',
\name{y}'')$ over $\mathop{\prod^{\bA^q}}\limits_{n<\omega}
\name{M}^n_{\bA^q}/\name{\cF}^q$ in the model $\mathop{{\prod}^{\bA^{q'}}}
\limits_{n<\omega} M^n_{\bA^q}/\name{\cF}^{q'}$ (i.e., the vocabulary and
the $\omega$ structures are from $\bV[G][\name{H}\cap {\mathbb C}_{\bA^q}]$,
the ultraproduct is taken in $\bV[G][\name{H}\cap {\mathbb C}_{\bA^{q'}}]$)
are equal.
\end{enumerate}
\underline{Then} the following conditions are equivalent.
\begin{enumerate}
\item[(A)] There is $r^0\in\App$ such that $q_\delta,q'\leq r^0$,
$r^0\restriction\delta\in G\cap(\App\restriction\delta)$, and  
\[\begin{array}{r}
p\forces_{{\mathbb C}_{\bA^{r^0}}}\mbox{`` }\mathop{\prod^{\bA^{r^0}}}
\limits_{n<\omega}\name{N}^1_n/\name{\cF}^{r^0}\models\ytheta[\name{x}'/
\name{\cF}^{r^0},\name{x}_\delta/\name{\cF}^{r^0}]\  \mbox{ and }\ \\
\mathop{\prod^{\bA^{r^0}}}\limits_{n<\omega}\name{N}_n^2/\name{\cF}^{r^0}
\models\neg\ytheta[\name{y}'/\name{\cF}^{r^0},\name{y}^*/\name{\cF}^{r^0}]
\mbox{ ''}. 
  \end{array}\]
\item[(B)] There is $r^1\in\App$ such that $q_\delta,q'\leq r^1$,
$r^1\restriction\delta\in G\cap(\App\restriction\delta)$ and  
\[\begin{array}{r}
p\forces_{{\mathbb C}_{\bA^{r^1}}}\mbox{`` }\mathop{\prod^{\bA^{r^1}}}
\limits_{n<\omega}\name{N}^1_n/\name{\cF}^{r^1}\models\ytheta[\name{x}''/
\name{\cF}^{r^1},\name{x}_\delta/\name{\cF}^{r^1}]\ \mbox{ and }\ \\
\mathop{\prod^{\bA^{r^1}}}\limits_{n<\omega}\name{N}_n^2/\name{\cF}^{r^1}
\models\neg\ytheta[\name{y}''/\name{\cF}^{r^1},\name{y}^*/\name{\cF}^{r^1}] 
\mbox{ ''}. 
  \end{array}\] 
\end{enumerate}
\end{lemma}

\begin{proof}
By symmetry it suffices to show that (A) implies (B). So suppose that $r^0$
is as in (A). We may also assume that $\bA^{r^0}=\bA^{q_\delta}\cup\bA^{q'}$ 
(just replace $q'$ by the stronger condition $r_\delta\restriction\delta$ if
needed). We want to define a respective condition $r^1$ so that $\bA^{r^1}=
\bA^{r^0}\stackrel{\rm def}{=}\bA$, and for this we need to extend
$\name{\cF}^{q_\delta}\cup\name{\cF}^{q'}$ to an ultrafilter containing the
set  
\[\{n\in\omega:\name{N}^1_n\models\ytheta [x''(n),\name{x}_\delta(n)]\}\cap  
\{n\in\omega:\name{N}^2_n\models\neg\ytheta [y''(n),\name{y}^*(n)]\}.\]
Suppose toward contradiction that this is impossible and thus we have a
condition $p'\in\coA$ stronger than $p$, a ${\mathbb C}_{
\bA^{q_\delta}}$--name $\name{\bf a}$ of a member of $\name{\cF}^{
q_\delta}$, and a ${\mathbb C}_{\bA^{q'}}$--name $\name{\bf b}$ of a member
of $\name{\cF}^{q'}$ such that $p'$ forces  
\[\mbox{`` }\name{\bf a}\cap\name{\bf b}\cap\{n\in\omega:\name{N}^1_n
\models\ytheta[\name{x}''(n),\name{x}_\delta(n)]\ \&\ \name{N}^2_n\models
\neg\ytheta[\name{y}''(n),\name{y}^*(n)]\}=\emptyset\mbox{ ''.}\]
Let $\bA_1=\bA^{q_\delta}$, $\bA_2=\bA^{q'}$, $\bA_0=\bA_1\cap\bA_2=\bA^q$,
and $p_i=p'\restriction \bA_i$ (for $i=0,1,2$). Let $H^0\subseteq{\mathbb
C}_{\bA_0}$ be generic over $\bV$ such that $p_0\in H^0$, and for
$n\in\omega$ let $\name{A}^*_n$ be the ${\mathbb C}_{\bA_0}$--name such that 
\[\begin{array}{ll}
\name{A}^*_n[H^0]=\big\{(u,v):&\mbox{there is a condition $p_1'\in {\mathbb
C}_{\bA_1}$ such that }\\
&p_1\leq p_1',\ p_1'\restriction \bA_0\in H^0\mbox{ and }p_1'\mbox{ forces
}\\
&\mbox{``}n\in\name{\bf a}\ \&\ \name{N}^1_n\models\ytheta[u,\name{x}_\delta 
(n)]\ \&\ \name{N}^2_n\models\neg\ytheta[v,\name{y}^*(n)]\mbox{'' }\big\}.
  \end{array}\]
Let $Q_{\name{\bar{A}}^*}$ be the predicate in $\tau_{\bA_0}$ corresponding
to the sequence $\name{\bar{A}}^*=\langle\name{A}^*_n:n<\omega\rangle$, see
\ref{generalia}(4). Note that   
\begin{enumerate}
\item[$(\otimes)'$] \qquad $p_2\forces_{{\mathbb C}_{\bA_2}}\mbox{`` }
\mathop{\prod^{\bA^{q'}}}\limits_{n<\omega} \name{M}^n_{\bA^q}/
\name{\cF}^{q'}\models Q_{\bA^*}(x',y')\mbox{ ''}$. 
\end{enumerate}
[Why? Assume not and let $H'\subseteq{\mathbb C}_{\bA}$ be a generic over
$\bV$ such that $H^0\subseteq H'$, $p'\in H'$, and we have (in $\bV^{H'}$)  
\[\name{\bf c}\stackrel{\rm def}{=}\{n:\neg \name{A}^*_n(x'(n),y'(n))\}\in
\name{\cF}^{q'}\subseteq\name{\cF}^{r^0}\quad\mbox{ and }\quad \name{\bf a}
\in \name{\cF}^{q_\delta}\subseteq\name{\cF}^{r^0}.\]
By our assumption
\[\name{\bf d}\stackrel{\rm def}{=}\{n: \name{N}^1_n\models \ytheta[
\name{x}'(n),\name{x}_\delta(n)]\ \&\ \name{N}^2_n\models\neg \ytheta[
\name{y}'(n),\name{y}^*(n)]\}\in\name{\cF}^{r^0},\]
so $\name{\bf a}\cap\name{\bf c}\cap\name{\bf d}\neq\emptyset$. Consequently
we may find $u,v$ and a condition $p^\otimes\in H'$ stronger than $p'$ such
that 
\[\begin{array}{ll}
p^\otimes\forces_{\coA}&\mbox{`` }x'(n)=u\ \&\ y'(n)=v\ \&\ \neg
A^*_n(u,v)\ \&\ n\in\name{\bf a}\ \&\\
&N^1_n\models\ytheta[u,\name{x}_\delta(n)]\ \&\ N^2_n\models\neg\ytheta[v,
\name{y}^*(n)]\mbox{ ''}.
  \end{array}\]
Then in particular $(u,v)\notin A^*_n[H^0]$, but also $p^\otimes\restriction
A_1$ witnesses $(u,v)\in A^*_n[H^0]$, a contradiction.]\\
Therefore, by assumption $(\beta)$, 
\begin{enumerate}
\item[$(\otimes)''$] \qquad $p_2\forces_{{\mathbb C}_{\bA_2}}\mbox{`` }
\mathop{\prod^{\bA^{q'}}}\limits_{n<\omega} M^n_{\bA^q}/\name{\cF}^{q'}
\models \name{A}^*(x'',y'')\mbox{ ''}$. 
\end{enumerate}
Thus we may choose $n,u,v$ and a condition $p_2'\in{\mathbb C}_{\bA^{q'}}$
such that $p_2\leq p_2'$, $p_2'\restriction \bA_0\in H^0$, and
\[p_2'\forces_{{\mathbb C}_{\bA_2}}\mbox{`` }n\in{\bf b}\ \&\ \name{x}''(n)
=u\ \&\ \name{y}''(n)=v\ \&\ \name{A}^*_n(u,v)\mbox{ ''}.\]
Since $(u,v)\in \name{A}^*_n[H^0]$, we may pick a condition $p'_1\in{\mathbb
C}_{\bA_1}$ stronger than $p_1$ and such that $p'_1\restriction \bA_0\in
H^0$ and 
\[p'_1\forces_{{\mathbb C}_{\bA_1}}\mbox{`` }n\in\name{\bf a}\ \&\
\name{N}^1_n\models\ytheta[u,\name{x}_\delta(n)]\ \&\ \name{N}^2_n\models
\neg\ytheta[v,\name{y}^*(n)]\mbox{ ''.}\] 
Then $p''\stackrel{\rm def}{=}p'\cup p'_1\cup p'_2\in{\mathbb C}_{\bA}$ is a 
condition stronger than $p'$  and 
\[p''\forces\mbox{`` }n\in\name{\bf a}\cap\name{\bf b}\ \&\ \name{N}^1_n
\models\ytheta[\name{x}''(n),\name{x}_\delta(n)]\ \&\ \name{N}^2_n\models
\neg\ytheta[\name{y}''(n),\name{y}^*(n)]\mbox{ '',}\] 
a contradiction.
\end{proof}

Let us go back to the proof of \ref{3.1D}. We define some ${\mathbb
C}_\delta$--names; recall $\name{H}^\delta\subseteq\coh\restriction\delta$
is generic over $\bV^*$, $\name{\cF}^\delta[\name{H}^\delta]=\bigcup\{
\name{\cF}^{r'}[\name{H}^\delta]:r'\in G_\delta\}$, and
\[\name{M}^*_\delta=\mathop{{\prod}^\delta}\limits_{n<\omega}
\name{M}^n_\delta/\name{\cF}^\delta,\quad\mbox{ and }\quad
\name{N}^\ell_\delta=\mathop{{\prod}^\delta}\limits_{n<\omega}
\name{N}^\ell_n/\name{\cF}^\delta\qquad\mbox{(for $\ell=1,2$)}.\]
Let
\[\begin{array}{l}
\name{Z}_\delta^1[\name{H}^\delta]=\Big\{(\name{x}/\name{\cF}^\delta,
\name{y}/\name{\cF}^\delta)\in\name{N}^1_\delta\times\name{N}^2_\delta:
\mbox{ there are a $\tau(t)$--formula $\ztheta\in\name{\Delta}$ and}\\
\quad\qquad\qquad \mbox{conditions $p\in{\mathbb C}_{\aleph_3}$ and
$r^0\in\App$ such that}\\
\quad\qquad\qquad p^*\leq p,\ p\restriction\delta\in H^\delta,\mbox{
$\name{x},\name{y}$ are ${\mathbb C}_{\bA^{r^0}\cap\delta}$--names, and}\\
\quad\qquad\qquad q_\delta\leq r^0,\ r^0\restriction\delta\in G\cap(\App
\restriction\delta),\mbox{ and}\\
\quad\qquad\qquad p\forces_{{\mathbb C}_{\bA^{r^0}}}\mbox{`` }
\mathop{\prod^{\bA^{r^0}}}\limits_{n<\omega}\name{N}^1_n/\name{\cF}^{r^0}
\models\ztheta[\name{x}/\name{\cF}^{r^0},\name{x}_\delta/\name{\cF}^{r^0}]\
\mbox{ and }\ \\ 
\qquad\qquad\qquad\qquad\mathop{\prod^{\bA^{r^0}}}\limits_{n<\omega}
\name{N}_n^2/\name{\cF}^{r^0}
\models\neg\ztheta[\name{y}/\name{\cF}^{r^0},\name{y}^*/\name{\cF}^{r^0}]
\mbox{ ''}\ \Big\},
\ \\
\name{Z}_\delta^0[\name{H}^\delta]=(\name{N}^1_\delta\times\name{N}^2_\delta)
\setminus\name{Z}_\delta^1.  
  \end{array}\]
Now, it follows from \ref{1.15Proposition} (and \ref{1.9Lemma}) that 
\begin{enumerate}
\item[$(\boxdot)_\delta$] in $\bV[G\cap(\App\restriction\delta)]
[\name{H}^\delta]$, {\bf if} the types realized by $(\name{x}'/
\name{\cF}^\delta,\name{y}'/\name{\cF}^\delta)$ and $(\name{x}''/
\name{\cF}^\delta,\name{y}''/\name{\cF}^\delta)$ over the model
$\mathop{\prod^{\bA^{q_\delta}\cap\delta}}\limits_{n<\omega}\name{M}^n_{
\bA^{q_\delta}\cap\delta}/\name{\cF}^{q_\delta\restriction\delta}$ in the
model $\mathop{\prod^\delta}\limits_{n<\omega}\name{M}^n_{\bA^{q_\delta}\cap 
\delta}/\name{\cF}^\delta$ are equal, {\bf then}
\[(\name{x}'/\name{\cF}^\delta,\name{y}'/\name{\cF}^\delta)\in
\name{Z}^0_\delta\quad\mbox{ if and only if }\quad (\name{x}''/
\name{\cF}^\delta,\name{y}''/\name{\cF}^\delta)\in \name{Z}^0_\delta.\]
\end{enumerate}
Now, most clauses of \ref{3.1D} should be clear; we say more on
(d)(vii,viii), for simplicity for $m=1$. 

We let $\name{R}_{\delta,1}=\name{Z}^0_\delta$, so clause
(d)(vii)$(\circledast)_2$ holds. 

Since $\name{F}$ is (an $\App*\coh$--name for) a $\name{\Delta}$--embedding
from $\prod\limits_{n<\omega}\name{N}^1_n/\name{\cF}$ onto $\prod\limits_{n<
\omega}\name{N}^2_n/\name{\cF}$, if $\name{x}/\name{\cF}^\delta\in
\name{N}^1_\delta$, then $\forces_{{\mathbb C}_\delta}\mbox{``
}(\name{x}/\name{\cF}^\delta,\name{F}(\name{x})/\name{\cF}^\delta)\in
\name{Z}^0_\delta\mbox{ ''}$. Hence clause (d)(viii)$(\oplus)_1$ holds.
\medskip

Thus the proof of \ref{3.1D} is completed. \hfill\hspace{0.2in}\vrule width
6pt height 6pt depth 0pt\vspace{0.1in} 
\medskip

\begin{conclusion}
\label{3.4}
In $\bV[G][H^{\aleph_3}]$, for each $\bm$, there is a stationary set
$S\subseteq\{\delta<\aleph_3:\cf(\delta)=\aleph_2\}$ and conditions
$q,q_\delta\in\App$ such that for each $\delta\in S$:  
\begin{itemize}
\item clauses (a)$_\delta$--(d)$_\delta$ of \ref{3.1C} are satisfied,
\item $q_\delta\in G$, $q_\delta\restriction \delta=q$,
$q_\delta,\name{y}_\delta$ as in \ref{3.1C},
\item the conclusion of \ref{3.1D} holds,
\item for every $\delta_1,\delta_2\in S$ there is a one-to-one order
preserving function $h:\bA^{q_{\delta_1}}\stackrel{\rm onto}{\longrightarrow}
\bA^{q_{\delta_2}}$ (so it is the identity on $\bA^q$) which maps
$\delta_1,\name{x}_{\delta_1},\name{F}(\name{x}_{\delta_1})=
\name{y}_{\delta_1}$ onto $\delta_2,\name{x}_{\delta_2},\name{F}(
\name{x}_{\delta_2})=\name{y}_{\delta_2}$,
\item in particular $p_\delta=p^*$ for $\delta\in S$, so the last clause in
\ref{3.1D} holds for $M_{\aleph_3}$.
\end{itemize}
\end{conclusion}

\begin{proof}
Straightforward.
\end{proof}

\section{Back to Model Theory}
In this section we present just enough to solve the problem on finite
fields. 

\begin{definition}
\label{4.1}
Let $M$ be a model. Assume $N_1=M^{\bar{\varphi}^1}$ , $N_2=
M^{\bar{\varphi}^2}$ are models of $t_0$ interpreted in $M$ by the sequences
$\bar{\varphi}^1,\bar{\varphi}^2$ of formulas with parameters from $M$, and
they have the same vocabulary $\tau^*=\tau(N_1)=\tau(N_2)$. Furthermore, let
$\Gamma$ be an invariant bigness notion in $M$ (over some set $A_0$ of
$<\kappa$ parameters, more exactly in $\cK_{(M,A_0)}$), and $\Delta\subseteq
L_{\omega,\omega}(\tau(N_1))$ and $\kappa>\aleph_0$ (for simplicity). 
\begin{enumerate}
\item We say that {\em $(N_1,N_2)$ is $(\kappa,\Gamma,\Delta)$--complicated
in $M$\/} if:\\
for every $\Delta$--embedding $F$ of $N_1$ into $N_2$, and for every
$\Gamma$--big type $\gp_0(x)$ inside $M$ of cardinality $<\kappa$,
\underline{there is} a $\Gamma$--big type $\gp_1(x)$ inside $M$ of
cardinality $<\kappa$ which includes $\gp_0(x)$ and such that, letting
$\tau(\gp)\subseteq\tau(M)$ consist of those mentioned in $\gp(x)$ (so
$|\tau(\gp)|<\kappa$) and $A\subseteq M$ be the set of parameters of $\gp_0$
(so $A_0\subseteq A$), we have 
\begin{enumerate}
\item[$(*)_{\gp(x)}$]  if 
\[\begin{array}{ll}
R_m=\{(\bar{a},\bar{b}):&\bar{a}\in {}^m(N_1),\ \bar{b}\in{}^m(N_2)\mbox{
and for some }\bar{c}\in{}^m(N_1)\mbox{ we have}\\
&{\rm tp}_{L_{\omega,\omega}(\tau(p))}(\bar{a}^\frown\bar{b},A,M)=
{\rm tp}_{L_{\omega,\omega}(\tau(p))}(\bar{c}^\frown F(\bar{c}),A,M)\ \}
  \end{array}\]
then we have the parallel of \ref{3.1D}(vii), so
\begin{enumerate}
\item[$(\oplus)_1$]  if $\bar{a}_1,\bar{a}_2$ are finite finite sequences of
the same length $m$ of members of $N^1_\delta$, and $\gp\cup\{\xtheta^{N_1}
(x,\bar{a}_1),\neg\xtheta^{N_1}(x,\bar{a}_2)\}$ is a $\Gamma$--big type
over $M_\delta$, and $\xtheta,\neg\xtheta\in\Delta$, \underline{then}
$(\bar{a}_1,F_\delta(\bar{a}_2))\notin R_m$. 
\item[$(\oplus)_2$] Above, we may replace $\xtheta,\neg\xtheta$ by any pair
$\xtheta_0,\xtheta_1$ of contradictory formulas from $\Delta$. 
\end{enumerate}
\end{enumerate}
\item In part (1):
\begin{enumerate}
\item[(i)] We do not mention $\Delta$ if it is the set of quantifier free
formulas (of $L_{\omega,\omega}(\tau(N_1))$).
\item[(ii)] We replace $\Gamma$ by $(t,\psi)$ if we mean ``for all bigness
notions of the form $\Gamma=\Gamma_{(t,\bar{\varphi},\psi)}$, where
$\bar{\varphi}$ is an interpretation of $t$ in $M$ with $<\kappa$ parameters
and $|t|<\kappa$, $\psi\in L_{\kappa,\omega}$'' (i.e., $\psi\in
L_{\mu^+,\omega}$ for some $\mu<\kappa$).
\item[(iii)] We omit $\Gamma$ if we mean ``for all $\Gamma$'s as in (ii)''.
\item[(iv)] We say {\em $M$ is $\kappa$--complicated\/} (or: {\em
$(\kappa,\Gamma,\Delta)$--complicated\/}) and omit $N_1,N_2$ \underline{if}
this holds for all $N_1,N_2$ as in our assumptions, but with
$|\tau(N_1)|<\kappa$.  
\end{enumerate}
\end{enumerate}
\end{definition}

\begin{remark}
We can add about the equivalence relation, implicit in \ref{3.1D}, see
\cite{Sh:F503}. 
\end{remark}

\begin{theorem}
\label{4.2}
\begin{enumerate}
\item Let $\bG$ be a full $(\aleph_3,\aleph_2)$--bigness guide (see
\ref{2.1A}; recall there is one by \ref{2.1B}). Assume that $G\subseteq
\App_\bG$ is generic over $\bV$ and $H\subseteq\coh$ is generic over
$\bV[G]$ and $\cF=\name{\cF}_{\aleph_3}[G][H]$, and let $\langle M_n=
M^n_{\aleph_3}:n<\omega\rangle$ be a sequence of models as above, that is
each with a countable universe being the set of natural numbers for
simplicity, all with the same vocabulary such that for every $k$  and
a sequence $\langle R_n:n<\omega\rangle$ with $R_n$ being a $k$--place
relation on $M_n$ there is a $k$-place predicate in the common vocabulary
satisfying $R^{M_n}=R_n$ for each $n$. \underline{Then} in $\bV[G][H]$ the
model $M=\prod\limits_{n<\omega} M^n_{\aleph_3}/\cF$ is
$\aleph_2$--complicated and $\aleph_2$--compact.
\item We can change the demands on $\bG$ accordingly to the version of
$\aleph_2$--complicated. 
\item If $N^1,N^2$ are models of $t^{\rm ind}_1$ interpreted in $M$, then
any isomorphism $\pi$ from $N^1$ onto $N^2$ is definable in $M$.
\item If $N^\ell=\prod\limits_{n<\omega} N^\ell_n/\cF$, each $N^\ell_n$ is  
countable, and $N^\ell$ is a model of $t^{\rm ind}_1$ (for $\ell=1,2$), and
$N^1,N^2$ are isomorphic, \underline{then} there are $A\in\cF$ and
isomorphisms $\pi_n$ from $N^1_n$ onto $N^2_n$ (for $n\in A$) such that
$\pi=\prod\limits_{n<\omega}\pi_n/\cF$.
\item  Above we replace : $N^\ell$ is a model of $t^{\rm ind}_1$"
by "for some  formula $\phi(x,y)$ in the vocabulary of   $N^1$
which is equal to that  of $N^2$, has the strong independence property
(in their common theory;
\footnote{ of course of the strong independence property
holds when we restrict ourselves to say a prediacte $P$ we get less,
but see \cite{Sh:F503}} ) 

\item If $N^\ell_n$ are finite fields (for $\ell=1,2$ and $n<\omega$), and
$\prod\limits_{n<\omega}N^1_n/\cF$ is isomorphic to $\prod\limits_{n<\omega}
N^2_n/\cF$, \underline{then} the set $\{n<\omega:N^1_n\simeq N^2_n\}$
belongs to $\cF$.
\end{enumerate}
\end{theorem}

\begin{proof}
(1)\quad By \ref{3.4}.
\medskip

\noindent (2)\quad The same proof.
\medskip

\noindent (3)\quad By \ref{4.3} below and \ref{1.3A}(2).
\medskip

\noindent (4)\quad Without loss of generality, the universe of $N^\ell_n$ is
$\alpha^\ell_n\leq\omega$. Now, for $\ell=1,2$, we can find
$P_\ell\in\tau_M$ such that $(P_\ell)^{M^n_{\aleph_3}}=|N^\ell_n|$ and for
$Q\in\tau(N^\ell_n)$ there is $Q_{[\ell]}\in\tau_M$ with $(Q_{[\ell]})^{
M^n_{\aleph_3}}=Q^{N^\ell_n}$. Therefore, $N^\ell=\prod\limits_{n<\omega}
N^\ell_n/\cF$ can be viewed as an interpretation in $M$ by
$\bar{\varphi}^\ell$. Now apply part (3) for $\Gamma=\Gamma_{(t^{\rm ind},
\bar{\varphi}^1,\psi^{\rm ind})}$. 
\medskip
\noindent (5)\quad  This follows by part (4), as the vacabulary
is finite, being a nn isomorphism is expresibleby a first order
sentence.
\noindent (6) \quad  This is a particular case of part (5).
By part (4) it suffice for infinite ultraproducts $N^\ell$ of finite fields
to find a formula $\xtheta(x,y)$ in the vocabulary of fields which has the
strong independence property.  First we deal with the case that the fileds
are of characteristic $p>2$. Consider the formula $\xtheta(x,y) $ saying
that $x+y$ has a square root in the field.

We relay a theorem of  Duret, \cite[p. 982, Lemma 10]{Du80}, chosing $p=2$,a
its hypothsis holds  as  the fieled contains all $p$-th root of the unit
(that is $1$  and $1$, The conlusion says that for  $n$  and any  pairwise
distinct elements $a_1,\dots, a_n, b_1,\dots, b_n$ of the field there is an
elemnt $c$ such that $a_m +c$ has a square root and $b_m +c$ does not have a
square root for $m=1,\dots, n$. So the formula $\xtheta(x,y)= (\exists z) (
z^2 = x+y)$ is as required. 

Of course, if the characteristic of the field is
$2$, then we naturally use the same theorem
but choosing $p=3$, so of course maybe the field
fail to have al $p$-th root of the unit,
however , as Duret does, we consider 
an algebraic extension of $N^\ell$ of order 3 
by adding a root of $X^`3-1$ hence all of them
getting a new filed $N^\ell_*$. Now the set of  
elements of $a^\ell_*N$
can be repreented as the set of triples of elements
of $N^\ell$, and the operations of $N^\ell_*$
 are definable in $N^\ell$,
So our problem is almost notational. E,g, we can note
that recalling $N^\ell == \prod_{n<\omega} N^\ell_n / \cF$
 then $N^\ell_*=\prod
N^\ell_{*,n} / \cF$  where $N^\ell_{*,n}$ is equal
 to $N^\ell_n$
if $N^\ell_n$ has three $3$-th roots of the unit
and an algebraic extension of $N^\ell_n$ of order three
which has this property otherwise.
Again the first order thoery of $N^\ell_*$ has the strong 
and for $N^1_*, N^\ell_2$ we get the desired conclusion
but any isomorphism from $N^1$  onto $N^2$
can be be  extended to an isomorphismfrom
$N^1_*$  onto $N^2_*$
and we can easily finish.

\end{proof}

\begin{proposition}
\label{4.3}
Assume that $M$ is a $\kappa$--complicated $\kappa$--compact model. Let  
$N_1,N_2$ be interpretations of $t^{\rm ind}_1$ in $M$. \underline{Then} for 
any isomorphism $\pi$ from $N_1$ onto $N_2$, the function $\pi\restriction
P^{N_1}$ is definable in $M$ by a first order formula (with parameters). 
\end{proposition}

\begin{proof}
Let $N_1=M^{\bar{\varphi}^\ell}$ (so $\bar{\varphi}^\ell$ has parameters in
$M$) and let $F$ be an isomorphism from $N_1$ onto $N_2$.

Let $\Gamma$ be the bigness notion $\Gamma_{(t^{\rm ind},\bar{\varphi}^1,
\psi^{\rm ind})}$ (so $\psi^{\rm ind}\in L_{\omega_1,\omega}$). Let
$A\subseteq M$, $|A|<\kappa$ and $\tau^*\subseteq\tau_M$, $|\tau^*|<\kappa$
be given by the definition of being $\kappa$--complicated (applied to
$F$). Without loss of generality, $A$ includes the parameters of
$\bar{\varphi}^1,\bar{\varphi}^2$ and  is closed under $F$ and $F^{-1}$, and
for every $n$ includes the finite set mentioned in \ref{bigind}(2).

Let $R_1$ be as in \ref{4.1}(1). Clearly, recalling Definition
\ref{bigind}(2), there are no distinct $a_1,a_2\in P^{N_1}\setminus A$ and
$b\in N_2$ such that $(a_1,b),(a_2,b)\in R_1$. Hence 
\[\{(b,a): (a,b)\in R_1\ \mbox{ and }\ a\in P^{N_1}\ \}\]
is the graph of a partial function from $P^{N_2}$ into $P^{N_1}$ which
includes the graph of $F^{-1}\restriction P^{N_2}$. Therefore,
$R_1\restriction (P^{N_1}\times P^{N_2})$ is the graph of $F\restriction
P^{N_1}$. But $R_1\restriction P^{N_1}$ is definable in $(M\restriction
\tau^*,c)_{c\in A}$ by a formula from $L_{\infty,\kappa}$, so also
$F\restriction P^{N_1}$ is, and thus if $N_1,N_2$ are models of $t^{\rm
ind}_1$ also $F$ is (by \ref{1.3A}). Applying \cite[1.9]{Sh:72} (or
\cite[Ch XI]{Sh:e}) 
we conclude that it is definable by a first order formula with parameters
from $M$, as required.   
\end{proof}

Similarly we can show the following.

\begin{proposition}
\label{4.4}
Assume that $\Gamma$ is a $(\aleph_2,\aleph_1)$--$(P,\xtheta)$--separative
bigness notion. Suppose that $N_1,N_2$ are interpretations of $t$ in $M$,
and $M$ is $\kappa$--compact $\kappa$--complicated (or just
$\kappa$--complicated for $\Gamma$), $\kappa>\aleph_0$. 
\begin{enumerate}
\item If $F$ is an isomorphism from $N_1$ onto $N_2$, \underline{then} 
\begin{enumerate}
\item[$(*)_1$] $F\restriction P^{N_1}$ is definable in $(M\restriction \tau^*,
c)_{c\in A}$ by a formula from $L_{\infty,\kappa}$, recalling $\tau\subseteq
\tau_M$, $|\tau|<\kappa$, $A\subseteq M$, $|A|<\kappa$.
\end{enumerate}
\item If $F$ is an embedding of $N_1$ into $N_2$, \underline{then} 
\begin{enumerate}
\item[$(*)_2$] there is a partial function $f$ from $P^{N_2}$ into $P^{N_1}$
which extends $F^{-1}$ and is definable in $(M\restriction\tau^*,c)_{c\in
A}$ by a formula from $L_{\infty,\kappa}$, where $\tau^*,A$ are as above. 
\end{enumerate}
\end{enumerate}
\end{proposition}

\begin{remark}
\begin{enumerate}
\item The proposition /ref{4.4}
 should be the beginning of an analysis of first order theories
$T$. For more in this direction see \cite{Sh:702}, \cite{Sh:F503}. 
\item As stated in the introduction, we may avoid the preliminary forcing
with $\App$ and construct the name $\cF$ in the ground model $\bV$, provided
$\bV$ is somewhat $\bL$--like. Assuming $\diamondsuit_{\{\delta\in\aleph_3:
\cf(\delta)=\omega_2\}}$ is enough, but we may also use the weaker principle
from \cite{HLSh:162} and \cite[Appendix]{Sh:405}.
\item We may vary the cardinals, e.g., we may replace $\aleph_2,\aleph_3$ by
$\kappa,\lambda$, respectively, provided $\lambda=2^\kappa$, $\kappa=
\kappa^{<\kappa}$ (so an approximation has size $<\kappa$).

So let us assume that
\[\theta=\theta^{<\theta}<\kappa=\kappa^{<\kappa}<\lambda=\kappa^+.\]
\begin{enumerate}
\item[(a)] For $\bA\subseteq \lambda$ let ${\mathbb C}(\bA)={\mathbb C}_\bA
=\{p: p$ is a partial function from $\Dom(p)\in [\bA]^{<\theta}$ to
${}^{\theta>} 2\ \}$ ordered by 
\[p_1\leq_{{\mathbb C}_\bA} p_2\quad\mbox{ iff }\quad \Dom(p_1)\subseteq
\Dom(p_2)\ \&\ (\forall\alpha\in\Dom(p_1))(p_1(\alpha)\trianglelefteq
p_2(\alpha).\] 
\item[(b)] We define $\App^-_\bG$ as the set of $q=(\bA^q,\name{\cF}^q)$ where
$A^q\in [\lambda]^{<\kappa}$ and $\name{\cF}^q$  is a ${\mathbb
C}_{\bA^q}$--name of a regular ultrafilter on $\theta$ such that for each
$\alpha<\lambda$, $\name{\cF}^q\cap{\mathcal P}(\theta)^{\bV^{{\mathbb
C}(\bA^q \cap \alpha)}}$ is a ${\mathbb C}_{\bA^q\cap \alpha}$--name.
\item[(c)] For $\alpha\in\bA\in [\lambda]^{<\kappa}$, $\name{x}_\alpha$ is
the ${\mathbb C}_\bA$--name $\bigcup\{p(\alpha):p\in\name{G}_{{\mathbb C}
(\bA)}$ of a member of ${}^\theta\theta$.
\item[(d)] We define $\name{N}^\varepsilon_\bA$ for $\varepsilon<\theta$,
$\bA\in [\lambda]^{<\kappa}$ as the following ${\mathbb C}_\bA$--name:

it is a model with universe $\theta$,

$\tau_{\name{N}^\varepsilon_\bA}=\{P_{\name{\bar{R}}}:
\name{\bar{R}}=\langle\name{R}_\varepsilon:\varepsilon<\theta$, for some $m$
each $\name{\bar{R}}_\varepsilon$ is an $m$--place relation $\}$,

$(P_{\name{\bar{R}}})^{\name{N}^\varepsilon_\bA}=\name{R}_\varepsilon$.

\noindent [So we may think of $\tau_{\name{N}^\varepsilon_\bA}$ to be an old
object whose members are indexed as $P_{\name{\bar{R}}}$, where each
$\name{R}_\varepsilon$ is a ${\mathbb C}_\bA$--name. Or we can consider
$\tau_{\name{N}^\varepsilon_\bA}$ to be a name and interpret it in
$\bV[G_{{\mathbb C}(\bA)}]$.
\end{enumerate}
\end{enumerate}
\end{remark}


\def\germ{\frak} \def\scr{\cal} \ifx\documentclass\undefinedcs
  \def\bf{\fam\bffam\tenbf}\def\rm{\fam0\tenrm}\fi 
  \def\defaultdefine#1#2{\expandafter\ifx\csname#1\endcsname\relax
  \expandafter\def\csname#1\endcsname{#2}\fi} \defaultdefine{Bbb}{\bf}
  \defaultdefine{frak}{\bf} \defaultdefine{mathfrak}{\frak}
  \defaultdefine{mathbb}{\bf} \defaultdefine{mathcal}{\cal}
  \defaultdefine{beth}{BETH}\defaultdefine{cal}{\bf} \def\bbfI{{\Bbb I}}
  \def\mbox{\hbox} \def\text{\hbox} \def\om{\omega} \def\Cal#1{{\bf #1}}
  \def\pcf{pcf} \defaultdefine{cf}{cf} \defaultdefine{reals}{{\Bbb R}}
  \defaultdefine{real}{{\Bbb R}} \def\restriction{{|}} \def\club{CLUB}
  \def\w{\omega} \def\exist{\exists} \def\se{{\germ se}} \def\bb{{\bf b}}
  \def\equivalence{\equiv} \let\lt< \let\gt> \def\cite#1{[#1]}
  \def\implies{\Rightarrow}

\shlhetal
\end{document}